\documentclass[12pt,draftcls,onecolumn]{IEEEtran}

\usepackage{makeidx}         
\usepackage{graphicx}        
\usepackage[pdftex]{color}
\usepackage{multicol}        
\usepackage[bottom]{footmisc}

\makeindex             
                       
\usepackage{amsfonts,amsmath,euscript}  

\newtheorem{theorem}{\bf Theorem}[section]

\newtheorem{definition}{\bf Definition}[section]
\newtheorem{lemma}{\bf Lemma}[section]

\newtheorem{problem}{\bf Problem}[section]
\newtheorem{proposition}{\bf Proposition}[section]
\newtheorem{remark}{\bf Remark}[section]
\newtheorem{example}{\bf Example}[section]

\newcommand{\bmat}{\left[ \begin{matrix}}
\newcommand{\emat}{\end{matrix} \right]}

\newcommand{\script}[1]{\EuScript{#1}}

\newcommand{\q}{\quad}

\newcommand{\Tr} {\mbox{\rm tr}}
\newcommand{\diag} {\mbox{\rm diag}}

\newcommand{\Var} {\mbox{\rm Var}\,}

\newcommand{\Circ}{\mathop{\rm Circ}}
\newcommand{\tp}{^{\top}}

\newcommand{\bea}{\begin{eqnarray}}
\newcommand{\eea}{\end{eqnarray}}

\newcommand{\bsea}{\begin{subeqnarray}}
\newcommand{\esea}{\end{subeqnarray}}

\newcommand{\Symmetric}{{\mathfrak S}}

\newcommand{\Ebb}{{\mathbb E}\,}

\newcommand{\Rbb}{\mathbb R}

\newcommand{\Zbb}{\mathbb Z}
\newcommand{\xb}{\mathbf  x}
\newcommand{\yb}{\mathbf  y}

\newcommand{\ab}{\mathbf a}
\newcommand{\bb}{\mathbf  b}

\newcommand{\db}{\mathbf  d}
\newcommand{\eb}{\mathbf  e}

\newcommand{\Cb}{\mathbf C}

\newcommand{\Fb}{\mathbf F}

\newcommand{\Ib}{\mathbf I}

\newcommand{\Mb}{\mathbf M}
\newcommand{\Nb}{\mathbf N}
\newcommand{\Qb}{\mathbf Q}

\newcommand{\Tb}{\mathbf T}
\newcommand{\Ub}{\mathbf U}
\newcommand{\Vb}{\mathbf V}

\newcommand{\Sigmab}{\boldsymbol{\Sigma}}

\newcommand{\Psib}{\boldsymbol{\Psi}}

\definecolor{Royalblue}{cmyk}{1,0.30,0.2,0.2}

\begin{document}

\title{ A Maximum Entropy solution of the Covariance Extension Problem for Reciprocal Processes}
\author{Francesca Carli, Augusto Ferrante, Michele Pavon, and Giorgio Picci\thanks{Work partially supported by the Italian Ministry for Education and Resarch (MIUR) under PRIN grant
``Identification and Adaptive Control of Industrial Systems".}
\thanks{F. Carli, A. Ferrante and G. Picci are  with the
Department of Information Engineering (DEI), University of Padova,
via Gradenigo 6/B, 35131
Padova, Italy. {\tt\small  carlifra@dei.unipd.it}, {\tt\small  augusto@dei.unipd.it},  {\tt\small  picci@dei.unipd.it}}
\thanks{M. Pavon is with   the
Department of Pure and Applied Mathematics, University of Padova,  {\tt\small
pavon@math.unipd.it}}}

\markboth{DRAFT}{Shell \MakeLowercase{\textit{et al.}}: Bare Demo of IEEEtran.cls for Journals}

\maketitle

\begin{abstract}Stationary {\em reciprocal processes}  defined on a  finite interval of the integer line  can  be seen as a special class of Markov random fields restricted to one dimension. Non stationary reciprocal processes  have been extensively   studied   in the past especially by Jamison, Krener, Levy and co-workers. The specialization of  the non-stationary theory to the stationary case, however,   does not seem  to have   been pursued  in sufficient depth in the literature.  Stationary reciprocal processes (and reciprocal stochastic models)  are potentially useful for describing  signals which naturally live in a finite region of the  time (or space) line.  Estimation or  identification of these models starting from   observed data seems still to be an open problem which can   lead to many interesting applications in signal and image processing. In this paper, we   discuss   a class of reciprocal processes which is the acausal analog of auto-regressive (AR) processes
 , familiar in control and signal processing. We show that  maximum likelihood identification of these processes leads to a covariance extension problem for block-circulant covariance matrices. This  generalizes  the famous  covariance band extension problem for  stationary  processes on the integer line. As in the usual stationary setting on the integer line,  the covariance extension problem turns out to be a basic conceptual and practical step in solving the identification problem. We show that the maximum entropy principle leads    to a complete solution of the problem. 
\end{abstract}

\section{Introduction}
  {\em Reciprocal processes}  have been introduced at the beginning of the last century \cite{Schrodinger-31,Bernstein-32, Schrodinger-32} even earlier than the idea of Markov process was formalized by Kolmogorov. The basic defining property is conditional independence given the values taken by the process at the boundary,  which resembles a widely  accepted definition of  Markov random fields. When the ``time'' parameter is  one dimensional, reciprocal processes can in fact  be seen as  Markov random fields restricted to one dimension. For this reason,   reciprocal processes are actually more general than Markov processes (a Markov process is reciprocal but not conversely). In fact,  these processes naturally live in a finite region of the time (or space) variable and specification of {\em boundary values} at the extremes of the interval is an essential part of their probabilistic description. In discrete-time they are naturally defined on a  finite interval of the integer l
 ine. Reciprocal processes  have been  extensively   studied   in the past notably by Jamison, Krener, Levy  and co-workers, see \cite{Jamison-70,Jamison-74,Jamison-75}, \cite{Krener-86,Krener-86b}, \cite{Levy-F-K-90}, \cite{Levy-F-02}, \cite{Frezza-90}. However  the specialization of  the non-stationary theory to the stationary case, except for  a few noticeable exceptions, e.g. \cite{Jamison-70}, \cite{Sand-94,Sand-96},  does not seem  to have   been pursued  in sufficient depth in the literature.  Stationary reciprocal processes (and reciprocal stochastic models)  are potentially useful for describing  signals which naturally live in a finite region of the  time or space line. They can be described by constant coefficient models which are a natural generalization of the Gauss-Markov state space models widely used in engineering and applied sciences. Estimation and  identification of these models starting from   observed data seems to be a completely open problem which can 
  lead to many interesting applications in signal and image processing. 
  
 In this paper,  after a general introduction to stationary processes defined on a finite interval (Section \ref{sec:CircCov}), we   discuss   a class of reciprocal processes  described by models which are the acausal analog of auto-regressive (AR) processes, familiar in control and signal processing (Section \ref{sec:circext}). In section \ref{sec:BilateralYW} we show that  maximum likelihood identification of these processes leads to a covariance extension problem for block-circulant covariance matrices. This  generalizes  the famous  covariance   extension problem for  stationary  processes on the integer line. As in the usual stationary setting on the integer line,  the covariance extension problem turns out to be a basic conceptual and practical step in solving the identification problem. The circulant covariance  extension problem    looks   similar to a classical    extension problems for positive block-Toeplitz matrices widely studied in the literature, \cite{Dym-G-81,Gohberg-G-K-94}, which belongs to the class of  {\em band extension} problems for positive matrices. All problems of this kind are solvable by factorization techniques.   However the banded algebra framework on which  this literature  relies     does not     apply to circulant matrices, see \cite{Carli-P-10}. Circulant band extension  appears to be a new kind of matrix extension problem.

In the present context, we are seeking a (reciprocal)  AR extension. One may speculate that this extension   should possess the  analog of the so-called ``maximum entropy'' property, which holds for stationary processes on the line.  In the literature, this property is usually  presented as a final  embellishment of the solution which is obtained by factorization techniques (typically computed via the Levinson-Whittle algorithm \cite{Levinson-47,Whittle-63}). In our case, where  there are no factorization techniques at hand, we resort to maximum entropy as the main tool at our disposal to attack the problem. In Sections \ref{MaxEntropy} and \ref{sec:VarAnal} we show that the maximum entropy principle indeed leads    to a complete solution of the problem. Finally in Section  \ref{sec:Reconciliation} we discuss  the relation  with the covariance selection results  in   Dempster's paper  \cite{Dempster-72}. 
  
Band extension problems for block-circulant matrices    of the type discussed in this paper  occur in particular  in applications to image modeling and simulation.  For reasons of space, we   do not provide details but rather refer the reader to the literature, see e.g. \cite{Chiuso-F-P-05,Chiuso-P-08} and \cite{Picci-C-08}.

\section{Stationary  processes on a finite interval}\label{sec:CircCov}

In this paper, we  work in the  wide-sense setting of second-order, zero-mean random variables. For the benefit of the reader, we recall here that a second order random   vector (or more generally process) is  just an equivalence class consisting of all zero-mean random vectors (or processes), each defined   on  some canonical probability space, say the space of their sample values,  that have the same covariance matrix, see e.g. \cite[Chap. X ]{Loeve-63}.  Hence, each second order random vector contains in particular a Gaussian element which may be taken as   the representative of the equivalence class, \cite[p. 74]{Doob-53}.  All statements of this paper do therefore apply to the particular case of Gaussian distributions. In  our setting, however,    explicit assumptions of Gaussianness will not be needed. We also recall that  there is  a basic correspondence,  established    by Kolmogorov in the early 1940's, between probabilistic concepts  depending only on second order m
 oments and geometric operations on certain  subspaces of  the Hilbert space of finite variance random variables,  see e.g. \cite[p. 636-637]{Doob-53} for historical remarks on this.  We assume henceforth that the reader is familiar with this correspondence.  

 Orthogonality of two random vectors will be understood as componentwise uncorrelation, i.e. $\xb \,\perp\, \yb$  means $\Ebb\xb\,\yb^{\top}=0$.  The symbol $\hat{\Ebb} [\,\cdot \mid \cdot\,]$   denotes  orthogonal projection (conditional expectation in the Gaussian case) onto the subspace spanned by a family of finite variance random variables listed  in the second argument.

A $m$-dimensional stochastic process on a finite interval $[\,1,\,N]$,   is just an ordered collection of  (zero-mean) random $m$-vectors  ${\yb}:= \{ \yb (k),\, k=1,2, \ldots,N \}$ which will be  written  as a column vector with $N$, $m$-dimensional components. We  say that $\yb$ is {\em stationary} if the covariances $\Ebb\yb (k)\yb (j)^{\top}$ depend  only on the difference of the arguments, namely
$$
\Ebb\yb (k)\yb (j)^{\top} = \Sigma_{k-j}\,, \qquad k,j= 1,\ldots,N,
$$
in which case the covariance matrix of $\yb$ has a symmetric block-Toeplitz structure; i.e.\footnote{Boldface capitals, e.g.  $\mathbf{I}_N$, $\Sigmab_{N}$, etc. denote block  matrices made of $N$ blocks,   each of dimension   $m\times m$.}
\begin{equation}\label{Toeplitz}
\Sigmab_{N} := \Ebb {\yb}{\yb}^{\top} = 
\bmat \Sigma_0& \Sigma_1^{\top}&\ldots &\Sigma_{N-1}^{\top}\\
          \Sigma_1 & \Sigma_0& \Sigma_1^{\top} & \ldots \\
          \vdots         &   \ddots       &   \ddots      & \ddots  \\
            \Sigma_{N-1} & \ldots      &   \Sigma_1& \Sigma_0  \emat   
\end{equation}
Processes ${\yb}$ which  have a   positive definite covariance $\Sigmab_{N}$ are called of   {\em full rank} (or {\em minimal}).  In this paper, we shall usually deal with full rank processes.

\begin{definition} A {\em block-circulant} matrix with $N$ blocks,  is a finite block-Toeplitz matrix whose
block-columns (or equivalently, block-rows)  are shifted cyclically. 
 
It looks like
$$
\Cb_{N}    = 
\bmat C_0           & C_{N-1}&\ldots &\ldots&C_1\\
        C_1  &         C_0& C_{N-1}    &\ldots  & \ldots \\
        \vdots     &          &       \ddots  &     & \vdots  \\
         \vdots     &          &        &     \ddots & C_{N-1}  \\
         C_{N-1} & C_{N-2}     &\ldots &   C_1 & C_0  \emat   \,.
$$
where $C_{k}\in \Rbb^{m\times m}$. A block-circulant matrix $\Cb_{N}$ is fully specified by its first block-column (or row).  It will be denoted by
\begin{equation}\label{CircC_comp}
\Cb_{N} = \Circ \{ C_0 ,  C_1 , \dots,C_{N-1}\}.
\end{equation}
\end{definition} 

For an introduction to circulant matrices, we refer the reader to the monograph \cite{Davis-79}. 
Block-circulant matrices of a fixed size form a real vector space which is actually an algebra with respect to the usual operations of sum and matrix multiplication. The invertible elements of this algebra form a group.\\
Consider now a stationary process $\tilde{\yb}$ on the integer line $\Zbb$, which is  periodic of period $T$, i.e. a process  satisfying   $\tilde{\yb}(k + n T): = \tilde{\yb}(k)$ (almost surely)  for all $n\, \in
\Zbb $.     We can think of $\tilde{\yb}$ as a process indexed on  the {\em discrete circle} group,  $\Zbb_{T} \equiv \{1,2,\ldots,T\}$ with arithmetics mod~$T$ \footnote{Whence  $T+\tau=\tau$ so that $T$ plays the role of the zero element.}. Clearly, its covariance function $\tilde \Sigma$  must also be periodic of period $T$,  namely,    $\tilde \Sigma_{k+T}=\tilde \Sigma_{k}$  for all $k \in
\Zbb $. Hence, we may also  see the covariance sequence as a function  on  the isomorphic  discrete group $\tilde{\Zbb}_{T} \equiv \,\{\,0,\,T-1\,\}$ with arithmetics mod~$T$.   But more must be true.

\begin{proposition}\label{prop:CircCov} 
A (second order) stochastic   process $\yb$ on $[\,1,\,T\,]$ is the restriction to the interval  $[\,1,\,T\,]$ of a wide-sense stationary periodic process $\tilde{\yb}$ of period $T$ defined on $\Zbb$,  
if and only if its covariance matrix $\Sigmab_T$ is symmetric block-circulant.
\end{proposition}\label{prop:symmCov}
 
\begin{IEEEproof} (only if) 
Let $k\in [\,1,\,T\,]$. By assumption there is  an $m$-dimensional stationary process $\tilde \yb$ on the integer line $\Zbb$, which is  periodic of period $T$,  satisfying   $\tilde{\yb}(k + n T): = {\yb}(k)$ (almost surely)  for arbitrary $n\, \in
\Zbb $.  By wide-sense stationarity, the covariance function of $\tilde \yb$ must depend only on the difference of the arguments, namely 
$$
\tilde{\Sigma}_{k,j}:= \Ebb\tilde{\yb}(k)\tilde{\yb}(j)^{\top}= \tilde{\Sigma}_{k-j}\,, \qquad k,j= 1,\ldots,T.
$$
Moreover, it is a well-known fact that, for any wide-sense  stationary process the following symmetry relation holds
\begin{equation}\label{eqn:symm}
\tilde{\Sigma}_{-\tau}=\tilde{\Sigma}_{\tau}^\top \qquad \forall \tau \in \Zbb\,,
\end{equation}
 that is the covariance matrix of $\tilde{\yb}$ has a symmetric block-Toeplitz structure. 
Now since $\tilde{\yb}$ is periodic of period $T$,  its covariance function must also be periodic of period $T$; i.e. $\tilde{\Sigma}_{k+nT}=\tilde{\Sigma}_k$ for arbitrary $k,n \in \Zbb$. 
Assume, just to fix the ideas, that $T$ is an even number and consider the midpoint $k=\frac{T}{2}$ of the interval $[1,\,T]$. 
The periodicity combined with the symmetry property \eqref{eqn:symm} yields that
\begin{equation}\label{eqn:symm_midpoint}
\tilde{\Sigma}_{\frac{T}{2} + \tau} =\tilde{\Sigma}_{\frac{T}{2} + \tau-T} = \tilde{\Sigma}_{\tau-\frac{T}{2}}=\tilde{\Sigma}_{\frac{T}{2}-\tau}^{\top} 
\qquad \forall \tau \in \Zbb
\end{equation}
and since \eqref{eqn:symm_midpoint} holds  for $\tau=0,1,\dots , \frac{T}{2}-1$, we can say that the   function $ \tilde{\Sigma}$ must be   {\em  symmetric with respect to the midpoint $\tau=\frac{T}{2}$ of the interval}. 
Hence, we can conclude that the 
covariance matrix of the process $\tilde{\yb}$ restricted to  $[\,1,\,T\,]$; that is the covariance $\Sigmab_T$ of  $\yb$, is a symmetric block-circulant matrix, i.e. it must have the following structure
$$ 
\Sigmab_T = 
\bmat \tilde{\Sigma}_0 &\tilde{\Sigma}_1^{\top} &\ldots        & \tilde{\Sigma}_{\tau}^{\top}&\ldots             & \tilde{\Sigma}_{\tau} & \ldots &       \tilde{\Sigma}_1\\
           \tilde{\Sigma}_1 & \tilde{\Sigma}_0         &\tilde{\Sigma}_1^{\top} &\ddots           & \tilde{\Sigma}_{\tau}^{\top}& \ldots       &    \ddots     &    \vdots  \\
	 \vdots   &          &       \ddots    	&	\ddots&  &  \ddots  &              &  \tilde{\Sigma}_{\tau}  \\
  \tilde{\Sigma}_{\tau}& \ldots &    \tilde{\Sigma}_1   & \tilde{\Sigma}_0 &  \tilde{\Sigma}_1 ^{\top}   &  \ldots &        \ddots & \\
                   \vdots & \tilde{\Sigma}_{\tau} &      \ldots      & & \tilde{\Sigma}_0 &          & \ldots &    \tilde{\Sigma}_{\tau} ^{\top} \\                                                  
  \tilde{\Sigma}_{\tau}^{\top}   &           &       \ddots    &        & &          &    				  & \vdots   \\
 \vdots   &  \ddots  &      &            \ddots       &       &	\ddots	&  \ddots & \tilde{\Sigma}_1^{\top}   \\
\tilde{\Sigma}_1^{\top} &\ldots  & \tilde{\Sigma}_{\tau}^{\top}&\ldots  &\tilde{\Sigma}_{\tau} & & \tilde{\Sigma}_1  &\tilde{\Sigma}_0 \emat 
$$
which we write
\begin{equation}\label{CircR}
\Sigmab_T = \Circ\{ \tilde{\Sigma}_0,\,\tilde{\Sigma}_1^{\top} ,\,\ldots ,\, \tilde{\Sigma}_{\tau}^{\top},  \ldots ,\,\tilde{\Sigma}_{\frac{T}{2}} ,\ldots ,\, \tilde{\Sigma}_{\tau},\, \ldots ,\,\tilde{\Sigma}_1\}\,.
\end{equation}
Similarly, if $T$ is odd, it must hold  that  
$\tilde{\Sigma}_{\frac{T+1}{2}+\tau}=\tilde{\Sigma}_{\frac{T-1}{2}-\tau}^{\top}$, $\tau=0,1, \dots, \frac{T-1}{2}-1$ 
and $\Sigmab_T$ can be written as 
$$
\Sigmab_T = \Circ\{ \tilde{\Sigma}_0,\,\tilde{\Sigma}_1^{\top} ,\,\ldots ,\, \tilde{\Sigma}_{\tau}^{\top},  \ldots ,\,\tilde{\Sigma}_{\frac{T-1}{2}}^{\top}, \, \tilde{\Sigma}_{\frac{T-1}{2}} ,\ldots ,\, \tilde{\Sigma}_{\tau},\, \ldots ,\,\tilde{\Sigma}_1\}\, ,
$$
which proves the first part of the statement.

(if) We want to prove that if $\yb$ is a process defined on a finite interval $[1,T]$ with a symmetric block-circulant covariance matrix $\Sigmab_T$, 
then it admits  a wide-sense stationary periodic extension, $\tilde{\yb}$, defined on $\Zbb$ of period $T$.

Let $\tilde{\yb}$ be the process obained by periodically extending the process $\yb$ to the whole interger line $\Zbb$ by setting   $\tilde{\yb}(k + n T): = {\yb}(k)$    for arbitrary $n\, \in \Zbb $ and 
let us denote by $\tilde{\Sigmab}$ its (infinite) covariance matrix. Since $\tilde{\Sigmab}$ is a covariance matrix, it must be positive semidefinite.  What we  need to show is that it is a symmetric block-Toeplitz matrix. By definition, $\tilde{\Sigmab}$ is the covariance matrix of the infinite column vector formed by stacking $\tilde{\yb}(0),\tilde{\yb}(1),\ldots,\tilde{\yb}(T),\ldots,\tilde{\yb}(n T),\ldots$ in that order, it is formed by subblocks which replicate $\Sigmab_T$ to produce a square matrix of infinite size. Since $\Sigmab_T$ is symmetric block-circulant, then $\tilde{\Sigmab}$ is, 
in particular, symmetric block-Toeplitz, which implies that $\tilde{\yb}$ is stationary. This concludes the proof. 
\end{IEEEproof}
\begin{remark}\label{rem:ext_symm_cov2} 
 The periodic  extension  to the whole line $\Zbb$ of deterministic signals originally given on a finite interval $[\,1,\,T]$ is  a common device in (deterministic) signal processing. This simple  periodic extension does however not preserve the  structure of a stationary random process  since the covariance of a periodically  extended process will not be 
stationary unless the  covariance function of the original process on $[\,1,\,T]$ was center-symmetric to start with. This counter-intuitive fact has to do with the quadratic dependence of the covariance of the process on its random variables. 
 
  Let for example $\yb$ be a scalar process on the finite interval $[1,4]$; i.e. let  $T=4$ and $m=1$.
Suppose $\yb$ has covariance matrix $\Sigmab_T = \text{Toepl}\left\{\sigma_0, \sigma_1, \sigma_2, \sigma_3\right\}$, the notation $\text{Toepl}\left\{a\right\}$ meaning that $\Sigmab_T$ is a symmetric Toeplitz matrix with first column given by the vector $a$.
 The upper-left $2T\times 2T$ corner the covariance of the periodic extension of $\yb$ is
 $$
\bmat
\sigma_0 & \sigma_1 & \sigma_2 & \sigma_3 & \sigma_0 & \sigma_1 & \sigma_2 & \sigma_3\\
\sigma_1 & \sigma_0 & \sigma_1 & \sigma_2 & \sigma_1 & \sigma_0 & \sigma_1 & \sigma_2\\
\sigma_2 & \sigma_1 & \sigma_0 & \sigma_1 & \sigma_2 & \sigma_1 & \sigma_0 & \sigma_1\\
\sigma_3 & \sigma_2 & \sigma_1 & \sigma_0 & \sigma_3 & \sigma_2 & \sigma_1 & \sigma_0\\
\sigma_0 & \sigma_1 & \sigma_2 & \sigma_3 & \sigma_0 & \sigma_1 & \sigma_2 & \sigma_3\\
\sigma_1 & \sigma_0 & \sigma_1 & \sigma_2 & \sigma_1 & \sigma_0 & \sigma_1 & \sigma_2\\
\sigma_2 & \sigma_1 & \sigma_0 & \sigma_1 & \sigma_2 & \sigma_1 & \sigma_0 & \sigma_1\\
\sigma_3 & \sigma_2 & \sigma_1 & \sigma_0 & \sigma_3 & \sigma_2 & \sigma_1 & \sigma_0\\
\emat.
$$
This matrix is clearly not Toeplitz unless $\sigma_3=\sigma_1$, in which case  $\Sigmab_T$ would be symmetric circulant. Hence the extended process $\tilde{\yb}$ is in general  not stationary. 
\end{remark}

\begin{remark}\label{rem:ext_symm_cov}   In many  applications to signal and image processing, the signals under study naturally live on a finite interval of the time (or space) variable and modeling them as functions defined on the whole line appears just as an artifice introduced in order to use the standard tools of (causal) time-invariant systems and harmonic  analysis on the line.   It may indeed be more logical to describe these data as    stationary processes $\yb$ defined on a finite interval $[1,T]$. The   covariance function,   say $\Sigmab_T$, of such a  process will    be a symmetric positive definite block-Toeplitz matrix which has in general no block-circulant structure. \\
It is   however always possible to extended the covariance function of $\yb$ to a larger interval so as to make it center-symmetric. This  can be achieved   by simply  letting $\Sigma_{T+\tau}:=  \Sigma_{T-1-\tau}^{\top}$ for $\tau=0, 1,\ldots ,T-1$.   In this way $\Sigmab_T$ is extended to a symmetric block-circulant matrix $\tilde \Sigmab_T$ of dimension   $(2T-1) \times (2T-1)$,  but  this operation   does not necessarily preserve positivity.  Positivity of a symmetric, block-circulant extension, however, can always be  guaranteed provided   the extension is done on a suitably   large interval. The details   on how to construct such an  extension are postponed to Section \ref{MaxEntropy}, see the proof of Theorem \ref{feas-maintheo}.  The original process ${\yb}$ can then be seen 
as the restriction to the interval $[1,\,T]$ of an extended process, say $\tilde{\yb}$, which lives on an interval $[1,\, N]$ of length $N \geq 2T-1$.  Since the extended covariance is, in any case, completely determined by the entries of the original covariance matrix $\Sigmab_T$, any statistical estimate thereof can be computed from the variables of the original process ${\yb}$   in the interval 
  $ [1,\,T]$ (or  from their sample values).   Hence, there is  no need to know what the random vectors 
$\{\tilde{\yb}(k)\,;\, k=T+1,\ldots,N\}$  look like.  
Indeed, as soon as we are given the covariance of the process ${\yb}$ defined on $[\,1,\,T]$, even if we may not ever see  
(sample values of) the ``external'' random vectors  $\{\tilde{\yb}(k)\,;\, k=T+1,\ldots, N\}$, we would in any case  have a completely determined second-order description (covariance function) of $\tilde{\yb}$.

In this sense, one can  think of any stationary process $\yb$ given on a finite interval $[1,\,T\,]$ 
as the restriction to $[1,T]$ of a wide-sense stationary {\em periodic} process,  $\tilde{\yb}$, of period $N\geq 2T -1$, defined on the whole integer line $\Zbb$. This process     naturally lives on the ``discrete   circle'' $\Zbb_{N}$. Hence dealing   in our future study   with the periodic  extension $\tilde{\yb}$, instead of the original process $\yb$,  will entail no loss of generality.  \hfill $\Box$
\end{remark}

\section{AR-type reciprocal processes}\label{sec:circext}

  In this section, we  describe a class of random  processes on a finite interval which are a natural generalization of the reciprocal processes introduced in \cite{Levy-F-K-90}, discussed in \cite{Levy-F-02} and, for the stationary case,  especially in  \cite{Sand-94,Sand-96}, see also \cite{Frezza-90}. In a sense, they are an acausal  ``symmetric''  generalization  of auto-regressive (AR) processes on the integer line.

Let $\yb$ be a   zero-mean  $m$-dimensional stationary process on $[1,\,N]$ and let $\Sigmab_N$ denote its $mN\times mN$  covariance matrix. We  assume that $\Sigmab_N$ is  a symmetric block-circulant matrix, so that $\yb$ may be seen as a process on the discrete circle $\Zbb_N$. 
In line with  what argued in Remark \ref{rem:ext_symm_cov}, we may, if we wish so, imagine that the matrix $\Sigmab_N$ was obtained by extending  a positive block-Toeplitz   matrix as \eqref{Toeplitz} to make it  symmetric block-circulant.  Then $[1,\,N]$ will have to be identified with an enlarged interval on  which $\yb$ is the periodic extension of some underlying stationary process. \\
Let $n$ be a natural number such that   $N>2n$. This inequality will be assumed to hold throughout.
We introduce the notation   $\yb_{[t-n,\,t\,)}$ for  the $nm$-dimensional random vector obtained by stacking $\yb(t-~n), \ldots, \yb(t-1)$ in that order. Similarly,   $\yb_{(t ,t+n\,]}$ is  the vector obtained by stacking $\yb(t+1), \ldots, \yb(t+n)$ in that order. Likewise, the vector  $\yb_{[t-n,\,t\,]}$ is obtained by appending  $\yb(t)$ as  last block to  $\yb_{[t-n,\,t\,)}$, etc.. The sums $t-k$  
and $t+k$   are to be understood modulo $N$. Consider a subinterval  $(t_1,\,t_2\,)\subset [1,\,N]$ where $(t_1,\,t_2\,):= \{t\mid t_1<t<t_2\}$ and $(t_1,\,t_2)^{c}$ denotes the complementary set in $[1,\,N]$. \\
Let $\script{A},\,\script{B},\,\script{C}$ be subspaces of zero mean second order random variables in a certain common ambient Hilbert space. Recall that $\script{A}$ and $\script{B}$   are said to be {\em conditionally orthogonal, given $\script{C}$} if 
\begin{equation} \label{condorthPerp}
\left(\ab -\hat{\Ebb}\,[\, \ab \mid \script{C}\,]\right) \,\perp\, \left(\bb -\hat{\Ebb}\,[\, \bb \mid\script{C}\,]\right) \,, \quad \forall\, \ab\,\in \script{A},\, \forall\, \bb\, \in \script{B}\,.
\end{equation}
Conditional orthogonality is the same as conditional uncorrelatedness (and hence conditional independence)  in the Gaussian case. Various equivalent forms of this condition are discussed in \cite{Lindquist-P-85}. When $\script{A},\,\script{B},\,\script{C}$ are  generated by finite dimensional random vectors,  condition \eqref{condorthPerp} can equivalently be rewritten in terms of the generating vectors, which we shall normally do in the following. The following definition does not require stationarity.

\begin{definition}\label{def:Recn}
A {\em   reciprocal process of order $n$}  on $ [1,\,N]$ is  characterized by the property that the random variables of the process in the   interval  $(t_1,\,t_2\,) $ are conditionally orthogonal  to the random variables in the exterior, $(t_1,\,t_2)^{c}$, given the $2n$ boundary values 
$\yb_{(t_1-n,\,t_1\,]}$ and $\yb_{[t_2 ,\,t_2+n\, )}$.
Equivalently,  it must hold that
\begin{equation}
\hat{\Ebb} [\,\yb_{(t_1,\,t_2)} \mid \yb(s),\, s\in (t_1,\,t_2)^{c}\, ]
= \hat{\Ebb} [\,\yb_{(t_1,\,t_2)}\mid \yb_{(t_1-n,\,t_1\, ]} \vee \yb_{[t_2 ,\,t_2+n\,)}\, ] \,, \quad t_1,\,t_2 \in\,[1,\,N] \,. 
\end{equation}
\end{definition}
In particular, we should have 
\begin{equation}\label{eqn:eq5pag8}
\hat{\Ebb} [\,\yb(t)\mid  \yb(s),\, s \neq t\, ]  = \hat{\Ebb} [\,\yb(t)\mid \yb_{[t-n,t\,)}\vee\yb_{(t ,t+n\,]} ] \,, \quad t\,\in\,\,[1,\,N], \,  
\end{equation}
where the  estimation error
\begin{equation}  \label{finconjn}
        \db(t)   : =  \yb(t)-\hat{\Ebb}[\,\yb(t)\mid  \yb(s),\, s \neq t\, ] = \yb(t) -\hat{\Ebb} [\,\yb(t)\mid \yb_{[t-n,t\,)}\vee\yb_{(t ,t+n\,]} ] \,, \quad t\,\in\,[1,\,N]\, 
\end{equation}
must clearly be  orthogonal to all random variables $\{ \yb(s),\, s \neq t\,\}$; i.e. 
\begin{equation}\label{OrthRecipn}
\Ebb\yb(t)\,\db(s)^{\top}\,= \Delta \, \delta_{ts}, \qquad t, \,s \,\in\,[1,\,N]\,,
 \end{equation}
where $\delta$ is the Kronecker function and $\Delta$ is a square matrix.  The actual meaning   of $\Delta$ will be clarified a few lines below. In the   spirit of Masani's definition \cite{Masani-60},  $\db$ is called the (unnormalized) {\em conjugate process}\footnote{Also called {\em double-sided innovation}.}  of $\yb$. Since $\db(t+k)$ is a linear combination of the components of the random vector $\yb_{[t+k-n,\,t+k+n]}$, it follows from \eqref{OrthRecipn} that both $\db(t+k)$ and  $\db(t-k)$ are orthogonal to $\db(t)$ as soon as $k >n$. Hence the process $\{\db(t)\}$  has   correlation   bandwidth~$n$; i.e.
\begin{equation}\label{MARecipn}
\Ebb \db(t+k)\,\db(t)^{\top}\,= \,0 \q \mbox{\rm for}\q   n < |k|  < N-n, \qquad k \,\in [0,\,N-1]\,. 
\end{equation}
It follows from    (\ref{finconjn}) that a  reciprocal process of order  $n$ on $ [1,\,N]$,  can always be described by a linear double-sided recursion  of the  form
\begin{equation}\label{Recipn}
\sum_{k=-n}^{n} F_k \,\yb(t-k)\,=\,\db(t)\,, \qquad t\in [1,\,N]
 \end{equation}
 where the $ F_k$'s are $m\times m$ matrices, in general dependent on $t$,  with $F_0=I_m$ 
 and $\db$ a    process of  correlation   bandwidth~$n$, orthogonal to $\yb$ in the sense of \eqref{OrthRecipn}. In fact, it follows from \eqref{OrthRecipn} that $\Ebb\db(t)\,\db(t)^{\top}\,= \Delta$  and hence $\Delta$ is the variance matrix of $\db(t)$,   symmetric and positive semidefinite.\\
 Equation \eqref{Recipn} requires the specification of boundary values, which will be described in Theorem \ref{thm:StatToAR(n)} below.
 
\begin{lemma}\label{lem:uniqeF}
{\em If $\yb$ is stationary,  the  matrices $\{F_k\}$ in the representation \eqref{Recipn} do not  depend on~$t$. If $\yb$ is full rank, they are uniquely determined by the covariance  lags  of the process up to order $2n$.}  
\end{lemma}		
\begin{IEEEproof} The $\{F_k(t)\}$'s are   determined by the orthogonality condition $\db(t)  \perp \yb_{[t-n,t\,)}\vee\yb_{(t ,t+n\,]}$, which can be expressed as 
\begin{align}\label{SolYW}
\bmat F_{-n}(t)& \ldots &F_{-1}(t) & F _{1}(t) & \ldots & F _{n}(t)\emat \,& \bmat \Sigmab_n & \Qb_n \\ \Qb_n^{\top} & \Sigmab_n \emat\, = \nonumber \\  \, = 
- \bmat  \Sigma_{n}^{\top} &\ldots \Sigma_{1}^{\top}& \Sigma_{1}& \ldots \Sigma_{n} \emat  &
\end{align}
where 
\begin{equation}\label{RQ}
 \Sigmab_n:=   \bmat \Sigma_{0} &  \Sigma_{1}& \ldots &  \Sigma_{n-1}\\
									 \Sigma_{1}^{{\top}} & \Sigma_{0} &\ldots & \\
									\ldots & \ldots& \ldots&  \Sigma_{1} \\
									 \Sigma_{n-1}^ {\top} & \ldots&  \Sigma_{1}^{\top}  &\Sigma_{0}\emat \,, \qquad
\Qb_n:= \bmat \Sigma_{n+1} &\Sigma_{n+2}& \ldots & \Sigma_{2n}\\
									\Sigma_{n} &\Sigma_{n+1}& \ldots & \Sigma_{2n-1}\\
									\ldots & \ldots& \ddots & \ldots\\
									\Sigma_2& \ldots & \Sigma_{n}& \Sigma_{n+1}\emat.
\end{equation}
Note that, because of stationarity, none of the covariance matrices depends on $t$. The determinant of the large block-matrix  in \eqref{SolYW} is  a principal minor of order   $2n$ of $\Sigmab_N$.  If $\yb$ is full rank, it must be nonzero and  the matrix must be invertible.  Therefore the  matrices $\{F_k\}$   do not  depend on~$t$ and are uniquely determined.
\end{IEEEproof}

For stationary reciprocal processes  on $\Zbb_N$, the boundary-values to be attached to the linear model \eqref{Recipn} are a straightforward consequence of the fact that  $\yb$  has a stationary periodic extension to the whole axis $\Zbb$. 
  \begin{theorem}\label{thm:StatToAR(n)}
{\em  A stationary reciprocal process, $\yb$, of order $n$ on $\Zbb_N$ satisfies a  linear,  constant-coefficients difference equation of the type \eqref{Recipn},  associated to the $2n$ cyclic boundary conditions:
\begin{equation}\label{CyclBC}
 \yb(k) = \yb(N+k)\,;  \qquad k=-n+1, \dots, n\,.
\end{equation}
The model  can be rewritten in matrix form as 
  \begin{equation}\label{MatRecipn}
 \Fb_{N}\, \yb = \db\,.
 \end{equation}
where $\Fb_{N}$ is the  $N$-block banded circulant matrix of bandwidth $n$, 
 \begin{equation}\label{FbandCirc}
 \Fb_{N} := \Circ\{ I, \,F _{1}, \, \ldots, \, F _{n}, \, 0 , \,\ldots 0, \, F _{-n}, \, \ldots, \, F _{-1}\}\,.
 \end{equation}
 If the process is full rank  this description is unique.}
 \end{theorem}
\begin{IEEEproof}
By definition 
 $$
 \hat{\Ebb} [\,\yb(1)\mid  \yb(s),\, s \neq 1\, ]  = \hat{\Ebb} [\,\yb(1)\mid \yb_{[1-n,\,1\,)}\vee\yb_{(1 ,\,1+n\,]} ] \,,  
$$
which is  a linear function of  $\yb_{[1-n,\,1\,)}\vee\yb_{(1 ,\,1+n\,]}$, whereby we can express $\yb(1)$ as
$$
\yb(1)=- \tilde{F}_-\, \yb_{(1 ,\,1+n\,]} - \tilde{F}_+ \,\yb_{[1-n,\,1\,)} +\db(1)
$$ 
 for some  coefficient matrices $\tilde{F}_-,\,\tilde{F}_+$. The process $\yb$ has a periodic extension of period $N$ and hence the missing initial boundary vector 
 $\yb_{[-n+1,\,1\,)}$ is actually the same as $ \yb_{[N-n+1,\,N\,]}$, so that
$$
\yb(1)=- \tilde{F}_-\, \yb_{(1 ,\,1+n\,]} - \tilde{F}_+ \,\yb_{[N-n+1,\,N\,]} +\db(1)\,.
$$ 
 By stationarity, the various $m\times m$ blocks in the  matrices $\tilde{F}$ must satisfy the same system of equations \eqref{SolYW} which was derived by imposing the orthogonality condition $\db(t)  \perp \yb_{[t-n,t\,)}\vee\yb_{(t ,t+n\,]}$,  for all times $t$.  
 Since the solution   is unique, it must hold that $\tilde{F}_{k}= F_{k}\,,\; k= \pm n,\,\ldots,\pm 1 $     where the $F_{k}$'s  are the same block matrices introduced  before  for \eqref{Delta}. Hence, we have
$$
 \yb(1)=- \sum_{k=-n}^{-1} F_{k}\,\yb(1-k) - \sum_{k=1}^{n} F_{k}\, \yb(N-k+1) +\db(1)
 $$
 which is the first ($t=1$) block equation in \eqref{Recipn} once   the first set of boundary conditions in \eqref{CyclBC} is used to replace the missing random variables $\yb_{[1-n,\,1\,)}$. Similar expressions can be derived for $\yb(2),\,\ldots,\yb(n)$ and for $\yb(N-n),\,\ldots,\yb(N)$. From this it readily follows that $\yb$ satisfies \eqref{MatRecipn} where $\Fb_N$ has the banded circulant structure \eqref{FbandCirc}. 
  \end{IEEEproof}
  
    Using the  notations $ F_-$ and $ F_+$ for $\bmat F_{-n}& \ldots &F_{-1}\emat$ and $\bmat F _{1} & \ldots & F _{n}\emat$ respectively, the error covariance $\Delta = \Var \{\db(t)\}$ can be expressed as
\begin{equation}\label {Delta}
\Delta = \Sigma_0 -\bmat F_-& F_+\emat \,  \bmat \Sigmab_n & \Qb_n \\ \Qb_n^{\top} & \Sigmab_n \emat^{-1}\,\bmat F_-& F_+\emat ^{\top}\,.
\end{equation}
 The following proposition is a simple generalization of analogous statements in \cite{Levy-F-K-90}, \cite{Sand-96} for $n=1$.
 \begin{proposition}\label{prop:Delta}
{\em A stationary reciprocal process $\yb$ is full rank if and only if the variance matrix $\Delta$ of the conjugate process is  positive definite.}   
\end{proposition}
\begin{IEEEproof} 
 (if) Suppose $\Delta >0$. Multiplying both members of \eqref{MatRecipn} from the right by $\yb^{\top}$ and taking expectations, in virtue of the orthogonality relation \eqref{OrthRecipn}, we get
\begin{equation}\label{CovMatRecipn}
\Fb_{N}\, \Sigmab_N =\Fb_{N}\,\Ebb\yb\yb^{\top} = \Ebb\db\yb^{\top}= \diag\{\Delta,\ldots,\Delta\}.
\end{equation}
Thus $\Delta >0$ implies that the square matrices $\Fb_{N}$ and $\Sigmab_N$ are invertible which, combined with the positive semidefiniteness of $\Sigmab_N$, implies $\Sigmab_N>0$.

 (only if) Suppose now that $\Delta$ is only positive semidefinite. This implies that there exists 
$ 0 \neq a \in \Rbb^m  \text{ s.t. } \Ebb a^{\top} \db(t)\db(t)^{\top}a = 0$, i.e.  s.t.  
$a^{\top} \db(t) = 0  \text{ a.s.}$. This means that the scalar components of $\db(t)$ are linearly dependent, which, by \eqref{Recipn}, implies 
that $\yb(t-n), \dots, \yb(t), \dots, \yb(t+n)$ are linearly dependent. Thus $\Sigmab_N$ 
must be  singular, which contradicts the assumption $\Sigmab_N>0$. 
\end{IEEEproof} 

Solving  \eqref {CovMatRecipn} we can express the inverse as 
\begin{equation}\label {RMinv}
\Mb_{N}:= \Sigmab_N^{-1} =\diag\{\Delta^{-1},\ldots,\Delta^{-1}\}\Fb_{N} 
\end{equation}
 so that  $\Mb_{N}$ is symmetric block-circulant  and positive definite, being the inverse of a matrix with the same properties. 
Furthermore,  $M_{k} :=  \Delta^{-1} F_{k}\,,\, k=-n,\ldots,n$ and $M_0= \Delta^{-1}$, must  form a {\em center-symmetric sequence} of bandwidth $n$; i.e.\footnote{That  is to say that  model \eqref{Recipn} is {\em self-adjoint}.}  
\begin{equation}\label{SymmF}
  M_{-k}= M_{k}^{\top}\,, \qquad k=1,\ldots,n\, .
\end{equation}
If we normalize the conjugate process by setting 
\begin{equation}\label {Normalizd}
\eb(t):=  \Delta^{-1}\db(t)
\end{equation}
so that $\Var \{ \eb(t)\}=  \Delta^{-1}$, the model \eqref{Recipn} can be rewritten
\begin{equation}\label{RecipnSymm}
\sum_{k=-n}^{n} M_k\,\yb(t-k)\,=\,\eb(t)\,, \qquad t\in \Zbb_{N}
 \end{equation}
for which the  orthogonality relation \eqref{OrthRecipn} is replaced by 
\begin{equation}\label{OrthRecipnSymm}
\Ebb\yb \,\eb ^{\top}\,=  \,\Ib_N\,.
 \end{equation}
We shall now show that  $\Mb_{N}$ 
is actually  the covariance matrix of the normalized conjugate process $\eb$. For, by the normalization \eqref{Normalizd}, our reciprocal process $\yb$ satisfies the linear equation 
\begin{equation}\label{NormRecipn}
\Mb_{N}\,\yb=\eb
\end{equation}
which  implicitly includes the cyclic boundary conditions \eqref {CyclBC}. Multiplying this from the right by $\eb^{\top}$ and taking expectations, we get
$ 
\Mb_{N}\,\Ebb\{\yb\eb^{\top}\}=\Ebb\{\eb\eb^{\top}\}
$ 
which, in force of \eqref{OrthRecipnSymm}, yields  
\begin{equation}
\Var\{\eb\}=\Mb_{N}
\end{equation}
as announced. We see that {\em the inverse of the covariance matrix of a  full rank stationary reciprocal process of order $n$, must be a banded block-circulant matrix of bandwidth $n$}.
 
This is in fact a fundamental characterization of  stationary reciprocal processes of order $n$. To prove it, we need to take up the (inverse) question of {\em well-posedness}, namely if  an autoregressive  model of the form \eqref{Recipn} associated to the proper cyclic boundary conditions, determines uniquely a process $\yb$ which is stationary and reciprocal of order $n$.    

To this end we may just as well examine the equivalent normalized model \eqref{NormRecipn}. 

 \begin{theorem}\label{thm:AR(n)toRec}
{\em Consider a  linear  model \eqref{NormRecipn}  where $\Mb_N$ is a symmetric positive-definite banded block-circulant  matrix of bandwidth $n$ and the process $\{\eb(t)\,;\; t\,\in \, \Zbb_N\}$  is  a stationary  process on $\Zbb_N$  with   covariance matrix $\Mb_N$.  

Then there is a unique full rank stationary reciprocal process $\yb$ of order $n$, solution of \eqref{NormRecipn}. This process  satisfies the orthogonality condition \eqref{OrthRecipnSymm} and $\eb$ is its normalized conjugate process.}
 \end{theorem}
 \begin{IEEEproof} 
 Pick a finitely correlated process $\eb$ with  covariance matrix $\Mb_{N}$ (we can construct such a, say Gaussian, process on a suitable probability space) and let $\yb$ be a solution of the equation \eqref{RecipnSymm} with boundary conditions \eqref{CyclBC}, equivalently a solution of  \eqref{NormRecipn}. Then, since  $\Mb_N$ is invertible, the process $\yb$ is uniquely defined on the interval $[1,\,N\,]$, i.e.  there is  a unique random vector,   $\yb$, solution of \eqref{NormRecipn}. Let $ \Sigmab_N$ be its covariance matrix. 
 We have,
 $\Sigmab_N:=\Ebb\left[\yb \yb^\top\right] = \Ebb\left[\Mb_N^{-1}\eb \eb^\top \Mb_N^{-\top}\right] = \Mb_N^{-1}$, so that $\Sigmab_N$   is a symmetric positive-definite  block-circulant  matrix and the process $\yb$ is stationary on $\Zbb_N$ (Proposition \ref{prop:CircCov}).
 
 By multiplying \eqref{NormRecipn} by $\eb^{\top}$ and taking expectations, we find
 $
 \Mb_N \Ebb\{\yb\eb^{\top}\} = \Mb_N\,,
 $
so that $\Ebb\{\yb\eb^{\top}\}= \Ib_N$, or equivalently
$
\Ebb\{\yb(t)\eb(s)^{\top}\}=  \Ib_m\,\delta_{ts}.
$
Therefore, the orthogonality \eqref{OrthRecipnSymm} holds on $\Zbb_N$.

Next, we need to show that $\yb$ is reciprocal of order $n$. To this end we shall generalize an argument of \cite{Sand-96}.
Let  $s<t$ be two points in $[1,\,N\,]$ , which for the moment we choose  such that $t-n> s+n$, which is always possible since by assumption $N>2n$. Expanding \eqref{RecipnSymm} and rearranging terms,  we can write
 \begin{align}\label{MN}
&\bmat M_0 &   M_1^{\top} &\ldots      & M_{n}^{\top}    &    0      &\ldots  &    0  && &0    \\
               M_1 & M_0 & M_1^{\top}  &\ddots & M_{n}^{\top}  & 0  &                    & & &0 \\
	           \vdots   &          &       \ddots    &                 &    	&	\ddots&       &&&\vdots \\
                                                 M_{n} & \ldots &    M_1   & M_0 &  M_1 ^{\top}  & \ldots&  M_{n}^{\top}&0&&0\\
                                                              0 & M_{n} &  &   \ldots           & M_0& \ldots &           &    \ddots  &&   \\
                  \vdots &     &   &   \ldots           &                            & \ddots &                                 &&&0  \\                                                              
                           0 &          &                     &   \ldots           &        & \ldots &                              &&& M_{n}^{\top} \\
          & &          &    				& &   &           &&&  \vdots  \\
       0 &&&   \ddots       &  & \ddots&		& M_1 &  M_0 & M_1^{\top}    \\
 0&&&     0  &\ldots  &0  &M_{n} & \ldots& M_1 &M_0\emat\,
\bmat \yb(t)\\ \yb(t+1)\\ \vdots \\ \yb(t+n )\\ \vdots \\  \yb(s-n)\\ \vdots \\ \yb(s-1)\\ \yb(s)\emat =   \notag\\
&   \bmat \eb( s)\\ \eb( s+1)\\ \vdots  \\ \eb(s+n )\\ \vdots \\  \eb(t-n )\\ \vdots \\ \eb(t-1)\\ \eb(t)\emat 
     - \bmat   M_{n}   & \ldots &      & M_1  & 0&\ldots &0 \\
                  0     &   M_{n}  & \ldots     &   M_2 & 0&\ldots &0  \\
	          0   &  \vdots  &    \ddots &  & 0&\ldots &0   \\
                   0   &     &  \ddots &      M_n& 0&\ldots &0\\
                  0& \ldots &    0  & \ldots & 0&\ldots &0 \\
                  \vdots& \vdots &    \vdots & \vdots & \vdots&\vdots &\vdots \\
                  0& \ldots &    0  & \ldots & 0&\ldots &0 \\
          0&\ldots &0 &         &  M_{n}^{\top}       &  0  &          \\
            0& \ldots &    0  & \ldots & \vdots& &0 \\
                 0&\ldots &0 &   &   M_2^{\top}  &        \ddots     &   0  \\
              0&\ldots &0&                   &  M_1^{\top}   &\ldots  & M_{n} ^{\top} \emat\, 
\bmat \yb(t-n)\\ \yb(t-n+1)\\ \vdots \\ \yb(t-1)\\ \yb(s+1)\\  \vdots \\ \yb(s+n-1)\\ \yb(s+n)\emat
\end{align}
which can be  compactly rewritten as
\begin{equation}\label{CompactEq}
\tilde{\Mb} \,\yb_{[t,\,s\,]}= \eb_{[t,\,s\,]} - 
\bmat \Nb & 0\\ 
0 & 0\\
0 &\Nb\tp\emat \bmat  \yb_{[t-n,\,t)} \\ \yb_{(\,s,\,s+n]} \emat
\end{equation}
with an obvious meaning of the symbols. Note that $\tilde{\Mb} $ is non-singular,  its determinant   being a principal minor of $\Mb_N$, and hence  nonzero;    while the two random vectors on the right hand side are uncorrelated since all scalar components of
$ \eb_{[t,\,s\,]}$ are orthogonal to the linear subspace spanned by (the scalar components  of) 
$\{\yb(\tau)\,;\, \tau \in [t,\,s\,]^{c}\}$  and hence are in particular orthogonal to the boundary condition  vectors $\yb_{(\,s,\,s+n]},\, \yb_{[t-n,\,t)} $. Solving \eqref{CompactEq} we can express    $\yb_{[t,\,s\,]}$ as a sum of two linear functions of  
$\eb_{[t,\,s\,]} $ and of $\yb_{(\,s,\,s+n]} \vee \yb_{[t-n,\,t)}$ so that the orthogonal projection onto the linear subspace spanned by (the scalar components  of) $\{\yb(\tau)\,;\, \tau \in [t,\,s\,]^{c}\}$ results in a linear function of (the scalar components of) $\yb_{[t-n,\,t)} \vee \yb_{(\,s,\,s+n]}$ alone. 
This proves the conditional orthogonality of   $\yb_{[t,\,s\,]}$ to the other random variables of the process, given the boundary values $\yb_{[t-n,\,t)}\,,\, \yb_{(\,s,\,s+n]}$.\\
 The argument     remains valid also when the non overlapping condition $t-n> s+n$ does not hold; i.e.  for an arbitrary  interval  $[t,\,s\,]$ of the discrete circle $\Zbb_N$. For, when  $[t-n,\,t)$ and $(\,s,\,s+n]$ overlap clearly we have  $[t,\,s\,]^{c} \subseteq [t-n,\,t)\cup (\,s,\,s+n]$ and hence all random variables in the subspace spanned by 
$\{\yb(\tau)\,;\, \tau \in [t,\,s\,]^{c}\}$ are contained in the subspace spanned by the boundary conditions, say $\script{C}:= \{\yb(\tau)\,;\, \tau \in [t-n,\,t)\cup (\,s,\,s+n]\}$. This means that $\hat{\Ebb} [\,\yb(\tau)\mid  \script{C}\, ] = \yb(\tau)$, or equivalently that
 $$
 \yb(\tau)-\hat{\Ebb} [\,\yb(\tau)\mid  \script{C}\, ] =0\,, \qquad \tau \in [t,\,s\,]^{c}
 $$
so that the second member in \eqref{condorthPerp} is zero and hence  the orthogonality condition trivially holds.
\end{IEEEproof}

  From this result, we obtain the following fundamental characterization of reciprocal processes on the discrete group $\Zbb_N$.

\begin{theorem}\label{InvCircBanded}
{\em A nonsingular $mN\times mN$-dimensional matrix $\Sigmab_N$ is the covariance matrix of a reciprocal process of order $n$
 on the discrete group $\Zbb_N$ if and only if its inverse is a positive-definite symmetric block-circulant matrix which is banded of bandwidth $n$.}
\end{theorem}

 Note that  the second order statistics of both $\yb$ and $\eb$ are encapsulated in the covariance matrix $\Mb_N$. In other words, the whole auto-regressive model of $\yb$ is defined in terms of  the matrix $\Mb_N$. 
 Note also that this result makes the stochastic realization problem for reciprocal processes of order $n$ conceptually trivial. In fact, given the covariance matrix $\Sigmab_N$ (the external description of the process), assuming   that it is in fact the covariance matrix of such a process, the model matrix $\Mb_N$ can be computed by simply inverting $\Sigmab_N$. This is the simplest answer one could hope for.
The solution  requires  however a preliminary criterion to check     whether   a (full rank) symmetric block-circulant covariance matrix has a banded  inverse.  There seems to be no  simple known answer  to this question.

Finally, to make contact with the literature, we note that a  full rank   reciprocal process  of order $n$ can always be represented as  a linear memoryless function of a   reciprocal process of order  1. This reciprocal process, however,   need not be of  full rank. To see that this is the case,   introduce the vectors
\begin{equation}
\yb_t^{+} := \bmat \yb(t) \\ \vdots \\  \yb(t+n-1)\emat\,, \quad \yb_t^{-} := \bmat \yb(t-n+1) \\ \vdots \\  \yb(t)\emat \,. \quad
\end{equation}
Letting $  \xb(t)^{\top}:=\bmat (\yb_t^{-})^{\top} &  ( \yb_t^{+})^{\top} \emat $, we find  the representation
\begin{align}
\xb (t) &= \bmat F_+ & 0 \\ 0 & 0\emat \xb(t-1) + \bmat 0 & 0 \\ 0 & F_-\emat \xb(t+1) + \tilde{\db}(t)\\
\yb(t) & =  \bmat 0& \ldots & 0 & 1 & 1 & 0 & \ldots & 0 \emat\, \xb(t)
\end{align}
where $ F_-$ and $ F_+$ are  the block-companion matrices
$$
F_+ :=\bmat 0& I & 0 &\ldots &0\\0&0 & I &\ldots &0\\  & \ldots & & &I \\ - F_n & \ldots & & & - F_1 \emat \qquad F_-  :=\bmat  - F_{-1} & \ldots & & & - F_{-n} \\ I & 0 &\ldots &  &0\\ 0&I &0 &\ldots &0\\  & \ldots & &I &0 \emat
$$
 and  $ \tilde{\db}(t) = \frac12\left[0 \; \dots \; 0 \; \db(t)^\top \; \db(t)^\top \; 0 \; \dots \; 0\right]^\top
$ has a singular covariance matrix. This model is in general non-minimal \cite{Sand-96}.

\section{Identification }\label{sec:BilateralYW}

Assume  that  $T$ independent   realizations of one period  of the  process  ${\yb}$ are available\footnote{For example,   a ``movie'' consisting of $T$ successive  images of the same texture.} and let us denote the string of   sample values by $\underline{y}:=\left(y^{(1)},	\dots ,y^{(T)}\right)$. We want to solve the following 

\begin{problem}\label{prob:IDP} Given the observations $\underline{y}$ of a reciprocal process ${\yb}$ of (known) order $n$, 
estimate the parameters $\{M_k\}$ of the underlying  reciprocal model $ \Mb_N \yb= \eb\,.$ 
\end{problem}

Note first that if we are given  $2n+1$ covariance data $\{\Sigma_k\,;\, k=0,1,\dots, 2n\}$, the identification of an  order $n$ reciprocal process can be carried  out by a linear algorithm,  namely by solving  the Yule-Walker-type  system of linear equations \eqref{SolYW}.

This procedure is however unsatisfactory since,  due to the symmetry \eqref{SymmF},  there are actually only $n+1$ unknown $M_k$ to be computed. Hence, one would expect only $n+1$ covariance lags to be needed, while the system \eqref{SolYW} requires solving also for the negative order coefficients. Moreover, in practice,  the $\Sigma_k$'s    will have to be  estimated from  observed data and estimates of covariances  with a large lag $k$ will unavoidably be more uncertain and have a larger variance.

 In an attempt to get asymptotically efficient estimates for the $M_k$'s, we  consider maximum likelihood estimation. To this end, we set up a Gaussian likelihood function (which does not require to assume that 
  $\yb$ has a Gaussian distribution, see \cite[p. 112]{Hannan-D-88}), which uses the density function
$$
p_{(M_0,\ldots,M_n)}(y) = \frac{1}{\sqrt{(2\pi)^{ mN}{\rm
det}\left(\Mb_N^{-1}\right)}}{\rm
exp}\left(-\frac{1}{2}y^{\top} \Mb_N y\right),
$$
where $y\in \Rbb^{mN}$.
Taking logarithms and neglecting terms which do not depend on the parameters, one can rewrite this expression as
\begin{align}
\log p_{(M_0,\ldots,M_n)}(y) = &-\frac12 \log\,\det\,\left(\Mb_N^{-1}\right) -\frac{1}{2}\,{\rm {\Tr}}\left\{ \Mb_N \,y y^{\top}\right\} \end{align}
Assuming   that the  $T$  sample measurements    are independent, the
 log-likelihood function, depending on the $n+1$ matrix parameters $\{M_k\,;\, k=0,1,\ldots,n\}$, can be written
\begin{equation}\label{MaxLik} 
L(M_0,\ldots,M_{n})   =  
   {\rm log\;} {\rm det}\left(\Mb_N\right)  - \sum_{k=0}^{n} {\rm  \Tr}\left\{M_k\,T_k\left(\underline{y}\right)\right\}  
\end{equation}
where each matrix-valued statistic  $T_k(\underline{y}) $ has the structure of  a sample  estimate   of the lag $k$ covariance of the process.
For example,  $T_0$ and $T_1$ are given by:
\begin{eqnarray*}
 T_0\left(\underline{y}\right) & = & \frac{1}{T}\sum_{t=1}^T
\left\{ \sum_{k =0}^{N -1} y^{(t)}(k)\left[y^{(t)}(k)\right]^\top \right\}     \\
 T_1\left(\underline{y}\right) &=&
\frac{2}{T} \sum_{t=1}^T\left\{ \sum_{k=1}^{N-1} y^{(t)}(k-1) \left[y^{(t)}(k)\right]^{\top}\right\}\\
& +&  \frac{2}{T}\sum_{t=1}^T y^{(t)}(N-1)\left[y^{(t)}(0)\right]^{\top}
\end{eqnarray*}
 From exponential class theory \cite{Barndorff-78}, we see that the $T_k$'s are (matrix-valued) sufficient statistics.   
Indeed, we have the well-known  characterization that the   (suitably normalized) statistics $T_0,\,T_1,\ldots ,T_n$ 
are Maximum Likelihood estimators of their expected values, namely 
\begin{eqnarray}\label{MLest_CovSampl_n}
\hat{ \Sigma}_0& :=&\frac{1}{N}T_0  = \;\mbox{\rm M.L. Estimator of}\;\Ebb {\yb}(k) {\yb}(k)^{\top} \nonumber\\
& \vdots & \\
\hat{ \Sigma}_n& : = &\frac{1}{N}T_{n}  = \;\mbox{\rm M.L. Estimator of}\; \Ebb {\yb}(k+n) {\yb}(k)^{\top} \,.\nonumber
\end{eqnarray}
Let us now consider the following matrix completion problem, which, form now on, will be referred to as the \emph{block-circulant band extension problem}. 

\begin{problem}[Block-Circulant Band Extension Problem]\label{prob:BCBEP}{\em Given   $n+1$ initial data $m\times m$  matrices  $\hat \Sigma_0,\,\dots, \, \hat \Sigma_n$, 
complete them with a sequence $\Sigma_{n+1},\Sigma_{n+2},\ldots, \Sigma_{N-1},$ in such a way to form a positive definite symmetric block-circulant matrix $\Sigmab_N$ with a 
{\bf block-circulant banded inverse of bandwidth $n$}.  } 
\end{problem}

Note that the model parameters $(M_0,\,M_1,\ldots, M_{n})$ are the nonzero blocks of the (banded) inverse of  the covariance  matrix $\Sigmab_N$ of the process (Theorem \ref{InvCircBanded}). The invariance principle for maximum likelihood estimators \cite{Zehna-66} leads then to  the following statement. 

\begin{theorem}\label{MainThm}
 {\em The maximum likelihood estimates of $(M_0,\,M_1,\ldots, M_{n})$ are the nonzero blocks of the banded inverse of  the matrix $\hat \Sigmab_N$ solving  the block-circulant band extension problem with initial data the $n+1$  covariance estimates \eqref {MLest_CovSampl_n}.}
\end{theorem}

 Hence, solving the original identification problem \ref{prob:IDP} has been shown to lead to the solution of a block-circulant band extension problem. Note, however, that the extension problem \ref{prob:BCBEP} is nonlinear and it is hard to see what is going on by elementary means.  Below we give   a scalar example. 

\begin{example}
Let  $m=1,\,N=8,\,n =2$ and assume we are  given the 
covariance estimates $\hat{\sigma}_0,\,\hat{\sigma}_1,\, \hat{\sigma}_2$, forming a positive definite Toeplitz matrix. 
The three unknown coefficients in the reciprocal model \eqref{RecipnSymm} of order $2$ are scalars, 
denoted $m_0,\,m_1,\, m_2$. Multiplying  \eqref{NormRecipn} from the right by $\yb^\top$, we get $\Mb_N \Sigmab_N= \Ib_N$, which leads to 
$$
\bmat m_0 &  m_1  &  m_2    &   0      &0  &    0  & m_2 &       m_1 \\
            m_1 &  m_0 &   m_1  &  m_2  & 0  &  0    &      0       &  m_2  \\
	  m_2 &  m_1 &  m_0 &   m_1  &   m_2  &  0  &  0    &      0   \\
             0 &  m_2 &  m_1 &  m_0 &   m_1  &   m_2  &  0  &  0       \\
             0  &    0 &  m_2 &  m_1 &  m_0 &   m_1  &   m_2  &  0     \\                 
              0 & 0  &    0 &  m_2 &  m_1 &  m_0 &   m_1  &   m_2    \\   
            m_2 & 0 & 0 &0  &  m_2 &m_1& m_0 &m_1\\
             m_1 & m_2 & 0 & 0 &0  &  m_2 &m_1& m_0 \emat
           \bmat \hat{\sigma}_0 \\ \hat{\sigma}_1\\ \hat{\sigma}_2 \\x_3 \\x_4\\ x_3\\ \hat{\sigma}_2\\\hat{\sigma}_1\emat= \bmat 1\\0\\ \vdots \\ \\0 \\0 \emat\,,
$$
where $x_3:=\sigma_3=\sigma_5$ and  $x_4:=\sigma_4$ are the unknown extended covariance lags. 
Rearranging and eliminating the last three redundant equations, one obtains
\begin{eqnarray*}
m_0\hat{\sigma}_0 + 2 m_1 \hat{\sigma}_1+ 2 m_2 \hat{\sigma}_2  &=& 1\\
m_0\hat{\sigma}_1 +   m_1 (\hat{\sigma}_{ 0} + \hat{\sigma}_2) +   m_2 (\hat{\sigma}_1 +x_3) &=& 0\\
m_0\hat{\sigma}_2 +   m_1 (\hat{\sigma}_1 + x_3) +   m_2 (\hat{\sigma}_0 +x_4) &=& 0\\
m_0x_3 +   m_1 (\hat{\sigma}_2 + x_4) +   m_2 (\hat{\sigma}_1 +x_3) &=& 0\\
m_0x_4+ 2 m_1 x_3 +  2 m_2 \hat{\sigma}_2  &=& 0
\end{eqnarray*}
which is a system of five quadratic equations in five unknowns whose solution   already looks  non-trivial. 
It may be checked that, under positivity of the  matrix  $\text{Toepl}\{\hat{\sigma}_0,\,\hat{\sigma}_1,\, \hat{\sigma}_2\}$, 
it has a unique positive definite solution (i.e. making $\Mb_N$ positive definite).
\end{example}

At  first sight the circulant band extension problem of Theorem \ref{prob:BCBEP} recalls the classical band extension problems 
for Toeplitz matrices studied in \cite{Dym-G-81,Gohberg-G-K-94}, which is solvable by factorization techniques.
However, the banded algebra framework on which these papers rely does not  apply here. 
The circulant band extension problem seems to be a new (and harder) extension problem.\\
General covariance extension problems   are discussed in an illuminating paper by A. P. Dempster, \cite{Dempster-72}. 
  Notice, however, that Dempster's procedures,  having been conceived to solve a general covariance extension problem, do  not exploit the circulant 
structure of the present setting and are computationally very intensive even for small scalar instances. 
A possible approximate approach to the circulant band extension problem was proposed in \cite{Chiuso-F-P-05}.  
This approach, based on a result of B. Levy \cite{Levy-92}, exploits the fact that for $N\to \infty$ the problem becomes one of band extension  for infinite positive definite symmetric  block-Toeplitz matrices, for which satisfactory algorithms exist. 
  For $N$ finite however,  this approximation may in some cases turn out to be poor. 
In the next section, we  propose a new approach to the circulant band extension problem.

\section{Maximum entropy on the discrete circle}\label{MaxEntropy}

Dempster's paper, which deals with general, unstructured covariance matrices, only considers  Gaussian distributions. He solves the following extension problem: Characterize, among all  covariance matrices sharing a given set of entries, the one corresponding to the (zero-mean) maximum entropy Gaussian distribution.
For our purposes, a key observation is  Statement (b) in \cite[p. 160]{Dempster-72}. In our setting, it reads as follows.

\begin{proposition}\label{prop:b}
{\em  Assume feasibility of the covariance extension problem. Among all covariance extensions of the data $\hat{\Sigma}_0\,\ldots,\hat{\Sigma}_n$, there exists a unique such an extension  whose inverse's  entries are zero in all the positions complementary to those where the elements of the covariance are assigned. This extension corresponds to the Gaussian distribution with  maximum  entropy.}
\end{proposition} 
This  principle of  entropy maximization  will   lead us to a  new convex optimization procedure   for computing the  band extension.

 We hasten to remark that in this paper we are {\em not} restricting ourselves to the case of Gaussian distributions.  We shall consider $\Sigmab_N$ to be the matrix variance of a Gaussian distribution only for the purpose of  {\em interpreting} the following optimization problem in the   light of Dempster's result. The far reaching implications of our maximum entropy principle for general probability distributions is provided in  Theorem \ref{GeneralEntropy} below.

\subsection*{Notations}
     Let ${\Ub}_N $ denote the  block-circulant ``shift" matrix with $N\times N$ blocks,
$$ 
{\Ub}_N =\bmat
0&I_m&0&\dots&0\\
0&0&I_m&\dots&0\\
\vdots&\vdots&
&\ddots&\vdots\\
       0&0&0&\dots&I_m\\
       I_m&0&0&\dots&0
       \emat,
$$
where $I_m$ denotes the $m\times m$ identity matrix. Clearly, 
${\Ub}_N \tp {\Ub}_N ={\Ub}_N {\Ub}_N \tp =I_{mN}$; i.e. ${\Ub}_N $ is orthogonal. Note that a matrix $C$ with $N\times N$ blocks is   block-circulant  if and only if it commutes with ${\Ub}_N $, namely if and only if it satisfies
\begin{equation}\label{charactcirc}
{\Ub}_N \tp C{\Ub}_N =C.
\end{equation}

Recall that the {\em differential entropy} $H(p)$ of a probability density function $p$ on $\Rbb^n$ is defined by
\begin{equation}\label{DiffEntropy}
H(p)=-\int_{\Rbb^n}\log (p(x))p(x)dx.
\end{equation}
  In  case of a zero-mean Gaussian distribution $p$ with covariance matrix $\Sigmab_N$, we get
\begin{equation}\label{gaussianentropy}
 H(p)=\frac{1}{2}\log(\det\Sigmab_N)+\frac{1}{2}n\left(1+\log(2\pi)\right).
\end{equation}
Let $\Symmetric_N$ denote the vector space of {\em symmetric} matrices with $N\times N$ square blocks of dimension $m\times m$. Let $\Tb_n \in \Symmetric_{n+1}$ denote the  Toeplitz matrix of {\em boundary data}:
\begin{equation}\label{ToeplInitData}
\Tb_n=\bmat 
\Sigma_{0}&\Sigma_{1}\tp&\ldots &\Sigma_{n}\tp \\
\Sigma_{1}& \ldots & & \ldots\\
\ldots& \ldots & & \ldots\\
\Sigma_{n}&\ldots & &\Sigma_{0} \emat
\end{equation}
and  let $E_n$ denote the $N\times (n+1)$ block matrix
$$
E_n=\bmat I_m&0&\ldots&0\\0&I_m&\ldots&0\\ 
					  0&0& \ldots&\ldots\\
				\ldots& &0&I_m\\
				0&0& \ldots&0  \emat.
$$
\subsection*{The Maximum Entropy problem on $\Zbb_N$} 
Consider the following Gaussian {\em maximum entropy} problem (MEP) on the discrete circle:
\begin{problem}\label{MaxEntPROBLEM}
\begin{eqnarray}
&&\min\left\{-\Tr\log \Sigmab_N \mid \Sigmab_N \in\Symmetric_{N},\; \Sigmab_N >0\right\}\\&&{\rm subject \;to:}\nonumber\\&&E_n^{\top} \Sigmab_N  E_n=\Tb_n,\label{c1}\\&&{\Ub}_N \tp \Sigmab_N {\Ub}_N =\Sigmab_N. \label{c2}
\end{eqnarray}
\end{problem}
Recalling that $\Tr\log\Sigmab_N=\log\det\Sigmab_N$ and (\ref{gaussianentropy}), we see that the above problem indeed amounts to finding the maximum entropy Gaussian distribution  with a block-circulant covariance, whose first $n+1$ blocks are precisely $\Sigma_{0}, \ldots,\Sigma_{n}$. The circulant structure  is equivalent to requiring this distribution to be stationary on the discrete circle $\Zbb_N$. We observe that in this problem we are minimizing a strictly convex function on the intersection of a convex cone (minus the zero matrix) with a linear manifold. Hence we are dealing with  a convex optimization problem.

Note  that we are   not imposing that the inverse of the solution $\Sigmab_N$ of Problem \ref{MaxEntPROBLEM} should have a banded structure.  We shall see that, whenever  solutions exist,  this   property  will be  {\em automatically  guaranteed.}
 
The first question to be addressed is feasibility of (MEP), namely the existence of a positive definite, symmetric matrix $\Sigmab_N$ satisfying (\ref{c1})-(\ref{c2}). Obviously, $\Tb_n$ positive definite is a necessary condition for the existence of such a $\Sigmab_N$. In general it turns out that, under such a necessary condition, feasibility holds for $N$ large enough.
The idea  is that for $N\to \infty$, Toeplitz matrices can be   approximated arbitrarily well   by circulants (\cite{Gray-02,Tyrthyshnikov-96}) and hence existence of a positive block-circulant extension can be derived from the existence of positive extensions for Toeplitz matrices.
\begin{theorem}\label{feas-maintheo}
{\em Given the sequence $\Sigma_i\in\Rbb^{m\times m}$, $i=0,1,\dots,n$, such that
\begin{equation}\label{feas-tnpd}
\Tb_n=\Tb_n\tp>0,
\end{equation}
there exists $\bar{N}$ such that for $N\geq\bar{N}$, the matrix $\Tb_n$ can be extended to an $N\times N$ block-circulant, positive-definite symmetric matrix $\Sigmab_N$.}
\end{theorem} 
\begin{IEEEproof}
A fundamental result in stochastic system theory  is the so-called maximum entropy covariance extension. It states   that, under condition (\ref{feas-tnpd}), there exists a rational positive real function $\Phi_+(z)=\frac{\Sigma_0}{2} +C(zI-A)^{-1}B$ such that
\begin{enumerate}
\item
$A$ has spectrum strictly inside the unit circle.
\item
$\Sigma_i=CA^{i-1}B$, $i=1,2,\dots, n.$
\item
The spectrum $\Phi(z):=\Phi_+(z)+\Phi_+^\ast(z)$ is coercive, i.e.\footnote{  Here, and in the following, ${\rm j}$ denotes the imaginary unit $\sqrt{-1}$.}
\begin{equation}\label{feas-coerc}
\exists \varepsilon>0 {\rm \ such\  that\ } \Phi({\rm e}^{{\rm j}\vartheta})>\varepsilon I,\ \forall \vartheta\in[0,2\pi).
\end{equation}
\end{enumerate}
 In fact $\Phi(z)$ has no zeros on the unit circle since it can be expressed in the form $\Phi(z)= L_n(z^{-1})^{-1} \Lambda_n  L_n(z )^{-\top}$ where $L_n(z^{-1})$ is the $n-th$ Levinson-Whittle matrix polynomial (also called $n-th$ matrix Szeg\"o polynomial) of the block Toeplitz matrix $\Tb_n$, and $\Lambda_n=\Lambda_n^{\top}>0$; see \cite{Whittle-63},  \cite{Delsarte-G-K-78} and \cite{Youla-K-78}.
 
Let $\Sigma_i:=CA^{i-1}B$, $i=n+1,n+2,\dots$, so that
$\Phi_+(z)=\frac{\Sigma_0}{2}+\sum_{i=1}^\infty \Sigma_iz^{-i}$,  and define
\begin{equation}
\Sigmab_N:=\left\{\begin{array}{ll} 
{\rm Circ} \left(\Sigma_0, \Sigma_1\tp, \Sigma_2\tp, \dots, \Sigma_{\frac{N-1}{2}}\tp,\Sigma_{\frac{N-1}{2}},\Sigma_{\frac{N-1}{2}-1},\dots \Sigma_1\right),& N {\rm\ odd}\\
{\rm Circ} \left(\Sigma_0, \Sigma_1\tp, \Sigma_2\tp, \dots, \Sigma_{\frac{N-2}{2}}\tp,\Sigma_{\frac{N}{2}}\tp+\Sigma_{\frac{N}{2}}, \Sigma_{\frac{N-2}{2}},\Sigma_{\frac{N-2}{2}-1},\dots \Sigma_1\right),& N {\rm\ even}
\end{array}
\right.
\end{equation}
We need  to show that there exists $\bar{N}$ such that $\Sigmab_N>0$ for $N\geq\bar{N}$.
To this aim, notice that $\Sigmab_N$ is, by definition, block-circulant so that, a similarity transformation induced by a unitary matrix $\Vb$ reduces  $\Sigmab_N$ to a block-diagonal matrix:
$$
\Vb^\ast \Sigmab_N \Vb=\Psib_N:=\diag\left(\Psi_0,\Psi_1,\dots,\Psi_{N-1}\right),
$$
where $\Vb$ is the Fourier block-matrix whose $k,l$-th block    is 
$$
V_{kl}=1/\sqrt{N} \exp\left[-{ {\rm j}}2\pi(k-1)(l-1)/N\right]I_m
$$ and $\Psi_\ell$ are the coefficients of the finite Fourier transform of the first block row of $\Sigmab_N$:
\begin{equation}
\Psi_\ell=\Sigma_0+{\rm e}^{{\rm j}\vartheta_\ell}\Sigma_1\tp+\left({\rm e}^{{\rm j}\vartheta_\ell}\right)^2\Sigma_2\tp+\dots +\left({\rm e}^{{\rm j}\vartheta_\ell}\right)^{N-2}\Sigma_2+\left({\rm e}^{{\rm j}\vartheta_\ell}\right)^{N-1}\Sigma_1,
\end{equation}
with $\vartheta_\ell :=- 2\pi\ell/N$,  see e.g. \cite[Sec. 3.4]{Tee-05}.
Clearly, $\left({\rm e}^{{\rm j}\vartheta_\ell}\right)^{N-i}=\left({\rm e}^{{\rm j}\vartheta_\ell}\right)^{-i}$ and 
hence
\begin{equation}\label{feas-fppsi}
\Psi_\ell=\Phi\left({\rm e}^{{\rm j}\vartheta_\ell}\right)-\left[\delta\Phi_N \left({\rm e}^{{\rm j}\vartheta_\ell}\right)+\delta\Phi_N^\ast \left({\rm e}^{{\rm j}\vartheta_\ell}\right)\right]
\end{equation}
where,
\begin{equation}
\delta\Phi_N(z):=\sum_{i=h+1}^\infty \Sigma_iz^{-i}=\sum_{i=h+1}^\infty CA^{i-1}Bz^{-i}
=z^{-h}CA^h(zI-A)^{-1}B,\quad
h:=\left\{\begin{array}{ll} 
\frac{N-1}{2},& N {\rm\ odd}\\
N/2, & N {\rm\ even}
\end{array}
\right.
\end{equation}

Since $A$ is a stability matrix, if $N$, and hence $h$, is large enough, 
$\delta\Phi_N \left({\rm e}^{{\rm j}\vartheta_\ell}\right)+\delta\Phi_N^\ast \left({\rm e}^{{\rm j}\vartheta_\ell}\right)$ is dominated by $\varepsilon I$, i.e.
there exists $\bar{N}$ such that
\begin{equation}\label{feas-condfp}
\delta\Phi_N \left({\rm e}^{{\rm j}\vartheta_\ell}\right)+\delta\Phi_N^\ast \left({\rm e}^{{\rm j}\vartheta_\ell}\right)< \varepsilon I, \quad \forall \vartheta_\ell,\quad \forall N\geq\bar{N}
\end{equation}
 so that it readily follows from (\ref{feas-coerc}) and (\ref{feas-fppsi}) that
 if $N\geq\bar{N}$,
$\Psi_\ell>0$ for all $\ell$.
\end{IEEEproof}

We observe that, given $\Tb_n$, the triple $A,B,C$ can be explicitly computed so that we can compute $\varepsilon$ and $\bar{N}$ for which (\ref{feas-condfp}) holds.
In other words, Theorem \ref{feas-maintheo} provides a sufficient condition that can be practically tested. Similar bounds, but valid only for the scalar case,   were derived in \cite{Dembo-M-S-89}.

\section{Variational analysis}\label{sec:VarAnal}
We shall  introduce a suitable set of ``Lagrange multipliers" for our constrained optimization problem. Consider the linear map
$ 
A:\Symmetric_{n+1}\times\Symmetric_N\rightarrow\Symmetric_N
$ 
defined by 
$$ A(\Lambda,\Theta)=E_n\Lambda E_n\tp +{\Ub}_N \Theta {\Ub}_N \tp -\Theta, \quad (\Lambda,\Theta)\in \Symmetric_{n+1}\times\Symmetric_N.
$$
and define the set 
$$
{\cal L}_+:=\{(\Lambda,\Theta)\in (\Symmetric_{n+1}\times\Symmetric_N)\mid (\Lambda,\Theta)\in(\ker(A))^\perp, \left(E_{n}\Lambda E_n^{\top} +{\Ub}_N \Theta {\Ub}_N \tp -\Theta\right) > 0\}.  
$$
Observe that ${\cal L}_+$ is an open, convex subset of $(\ker(A))^\perp$. For each $(\Lambda,\Theta)\in {\cal L}_+$, we consider the unconstrained minimization of the {\em Lagrangian function}
\begin{eqnarray}\nonumber 
L(\Sigmab_N,\Lambda,\Theta)&:=&-\Tr\log\Sigmab_N+\Tr\left(\Lambda\left(E_n^{\top}\Sigmab_N E_n-\Tb_n\right)\right)+\Tr\left(\Theta\left({\Ub}_N \tp \Sigmab_N {\Ub}_N -\Sigmab_N\right)\right)\\&=&-\Tr\log\Sigmab_N+\Tr\left(E_n\Lambda E_n\tp \Sigmab_N\right)-\Tr\left(\Lambda\Tb_n\right)+\Tr\left({\Ub}_N \Theta {\Ub}_N \tp \Sigmab_N \right) \nonumber \\
&&-\Tr\left(\Theta \Sigmab_N \right)\nonumber
\end{eqnarray}
over ${\mathfrak S}_{N,+}:=\{\Sigmab_N \in\Symmetric_N,\;\Sigmab_N>0\}$.  For $\delta\Sigmab_N\in\Symmetric_N$, we get 
$$
\delta L(\Sigmab_N,\Lambda,\Theta;\delta\Sigmab_N)=-\Tr\left(\Sigmab_N^{-1}\delta\Sigmab_N\right)+\Tr\left(E_n\Lambda E_n\tp \delta \Sigmab_N\right)+\Tr\left(\left({\Ub}_N \Theta {\Ub}_N \tp -\Theta\right)\delta\Sigmab_N\right).
$$
We conclude that $\delta L(\Sigmab_N,\Lambda,\Theta;\delta\Sigmab_N)=0,\; \forall \delta\Sigmab_N\in\Symmetric_N$ if and only if
$$
\Sigmab_N^{-1}=E_n\Lambda E_n\tp +{\Ub}_N \Theta {\Ub}_N \tp -\Theta.
$$
Thus, for each fixed pair $(\Lambda,\Theta)\in {\cal L}_+$, the unique $\Sigmab_N^o$ minimizing the Lagrangian is given by
\begin{equation}
\Sigmab_N^o=\left(E_n\Lambda E_n\tp +{\Ub}_N \Theta {\Ub}_N \tp -\Theta\right)^{-1}.
\end{equation}
Consider next $L(\Sigmab_N^o,\Lambda,\Theta)$. We get
\begin{eqnarray}
\nonumber L(\Sigmab_N^o,\Lambda,\Theta)=-\Tr\log\left(\left(E_{n}\Lambda E_n^{\top} +{\Ub}_N \Theta {\Ub}_N \tp -\Theta\right)^{-1}\right)\\+\Tr\left[\left(E_{n}\Lambda E_n^{\top} +{\Ub}_N \Theta {\Ub}_N \tp -\Theta\right)\left(E_{n}\Lambda E_n^{\top} +{\Ub}_N \Theta {\Ub}_N \tp -\Theta\right)^{-1}\right]-\Tr(\Lambda\Tb_n)\\
=\Tr\log\left(E_{n}\Lambda E_n^{\top} +{\Ub}_N \Theta {\Ub}_N \tp -\Theta\right)+\Tr I_{mN}-\Tr\left(\Lambda\Tb_n\right).\nonumber
\end{eqnarray}
This is a strictly concave function on ${\cal L}_+$ whose maximization is the {\em dual problem} of (MEP).
We can equivalently consider the convex problem
\begin{equation}\min\left\{J(\Lambda,\Theta),(\Lambda,\Theta)\in{\cal L}_+\right\},
\end{equation}
where $J$ (henceforth called dual function) is given by
\begin{equation}J(\Lambda,\Theta)=\Tr\left(\Lambda\Tb_n\right)-\Tr\log\left(E_n\Lambda E_n\tp +{\Ub}_N \Theta {\Ub}_N \tp -\Theta\right).
\end{equation}

\subsection*{Existence for the dual problem}
The minimization of the strictly convex function $J(\Lambda,\Theta)$ on the convex set ${\cal L}_+$ is a  challenging problem as ${\cal L}_+$ is an {\em open} and {\em unbounded} subset of $(\ker(A))^\perp$. Nevertheless,  the following existence result in the Byrnes-Lindquist spirit, \cite{Georgiou-L-03},  \cite{Byrnes-L-07}, \cite{Ferrante-P-R-07} can be established.
\begin{theorem}\label{existencedual}
{\em The function $J$ admits a unique minimum point $(\bar{\Lambda},\bar{\Theta})$ in ${\cal L}_+$.}
\end{theorem}
In order to prove this theorem, we need first to derive a number of auxiliary results.
Let $\mathfrak{C}_N$ denote the vector subspace of block-circulant matrices in $\Symmetric_N$. We proceed to characterize the orthogonal complement of $\mathfrak{C}_N$ in $\Symmetric_N$.
\begin{lemma}\label{perpcirc}
{\em Let $M\in\Symmetric_N$. Then $M\in(\mathfrak{C}_N)^\perp$ if and only if it can be expressed as
\begin{equation}
M={\Ub}_N N{\Ub}_N \tp -N
\end{equation}
for some $N\in\Symmetric_N$.}
\end{lemma}
\begin{IEEEproof}
By (\ref{charactcirc}), $\mathfrak{C}_N$ is the kernel of the linear map from $\Symmetric_N$ to $\Symmetric_N$ given by $M\mapsto {\Ub}_N \tp M{\Ub}_N -M$. Hence, its orthogonal complement is the range of the adjoint map. Since
$$\Tr\left(({\Ub}_N \tp M{\Ub}_N -M)N\right)=\langle {\Ub}_N \tp M{\Ub}_N -M,N\rangle=\langle M,{\Ub}_N N{\Ub}_N \tp -N\rangle,
$$
the conclusion follows.
\end{IEEEproof}

Next we show that, as expected, feasibility of the primal problem (MEP) implies that the dual function $J$ is bounded below.
\begin{lemma}\label{boundedbelow}
{\em Assume that there exists $\bar{\Sigmab}_N\in{\mathfrak S}_{N,+}$ satisfying (\ref{c1})-(\ref{c2}). Then, for any pair $(\Lambda,\Theta)\in{\cal L}_+$, we have}
\begin{equation}J(\Lambda,\Theta)\ge mN+\Tr\log\bar{\Sigmab}_N. 
\end{equation}
\end{lemma}
\begin{IEEEproof} By (\ref{c1}), $\Tr(\Lambda\Tb_n)=\Tr(\Lambda E_n\tp \bar{\Sigmab}_NE_n)=\Tr(E_{n}\Lambda E_n^{\top} \bar{\Sigmab}_N)$. Using this fact and Lemma \ref{perpcirc}, we can now rewrite the dual function $J$ as follows
\begin{eqnarray}\nonumber
J(\Lambda,\Theta)&=&\Tr\left(\Lambda\Tb_n\right)-\Tr\log\left(E_{n}\Lambda E_n^{\top} +{\Ub}_N \Theta {\Ub}_N \tp -\Theta\right)\\&=&\Tr\left[\left(E_{n}\Lambda E_n^{\top} +{\Ub}_N \Theta {\Ub}_N \tp -\Theta\right)\bar{\Sigmab_N}\right]-\Tr\log\left(E_{n}\Lambda E_n^{\top} +{\Ub}_N \Theta {\Ub}_N \tp -\Theta\right).\nonumber
\end{eqnarray}
Define $M(\Lambda,\Theta)=\left(E_{n}\Lambda E_n^{\top} +{\Ub}_N \Theta {\Ub}_N \tp -\Theta\right)$ which is positive definite for $(\Lambda,\Theta)$ in ${\cal L}_+$. Then
$$J(\Lambda,\Theta)=\Tr\left(M(\Lambda,\Theta)\bar{\Sigmab}_N\right)-\Tr\log M(\Lambda,\Theta).
$$
As a function of $M$, this is a strictly convex function on ${\mathfrak S}_{N,+}$, whose unique minimum occurs at $M=\bar{\Sigmab}_N^{-1}$ where the minimum value is $\Tr(I_{mN})+\Tr\log\bar{\Sigmab}_N$.
\end{IEEEproof}

\begin{lemma}\label{atinfinity}
{\em Let $(\Lambda_k,\Theta_k), n\ge 1$ be a sequence of pairs in ${\cal L}_+$ such that $\|(\Lambda_k,\Theta_k)\|\rightarrow\infty$. Then also $\|A\left(\Lambda_k,\Theta_k\right)\|\rightarrow\infty$. It then follows that $\|(\Lambda_k,\Theta_k)\|\rightarrow\infty$ implies that $J(\Lambda_k,\Theta_k)\rightarrow\infty$.}
\end{lemma}
\begin{IEEEproof}
Notice that $A$ is a linear operator between finite-dimensional linear spaces.  
Denote by $\sigma_m$ the smallest singular value of the restriction of  $A$ to $(\ker A)^\perp$ (the orthogonal complement of $\ker A$). Clearly, $\sigma_m>0$, so that,
since each element of the sequence $(\Lambda_k,\Theta_k)$ is in $(\ker A)^\perp$, $\|A\left(\Lambda_k,\Theta_k\right)\|\geq \sigma_m \|(\Lambda_k,\Theta_k)\|\rightarrow\infty$.

Assume now that $\|A\left(\Lambda_k,\Theta_k\right)\|=\|\left(E_n\Lambda_k E_n\tp +{\Ub}_N \Theta_k {\Ub}_N \tp -\Theta_k\right)\|\rightarrow\infty$. Since these are all positive definite matrices and all matrix norms are equivalent, it follows that 
  $$
  \Tr\left(E_{n}\Lambda E_n^{\top} +{\Ub}_N \Theta {\Ub}_N \tp -\Theta\right)\rightarrow\infty.
  $$
 As a consequence, $\Tr\left(\left(E_{n}\Lambda E_n^{\top} +{\Ub}_N \Theta {\Ub}_N \tp -\Theta\right)\bar{\Sigmab}_N\right)\rightarrow\infty$ and, finally, $J(\Lambda_k,\Theta_k)\rightarrow\infty$ .
\end{IEEEproof}

We show next that the dual function tends to infinity also when approaching the boundary of ${\cal L}_+$, namely
\begin{eqnarray}\nonumber\partial{\cal L}_+:=\{(\Lambda,\Theta)\in (\Symmetric_{n+1}\times\Symmetric_N)| (\Lambda,\Theta)\in(\ker(A))^\perp, \left(E_{n}\Lambda E_n^{\top} +{\Ub}_N \Theta {\Ub}_N \tp -\Theta\right) \ge 0,\\ \det\left(E_{n}\Lambda E_n^{\top} +{\Ub}_N \Theta {\Ub}_N \tp -\Theta\right)=0\}.  \nonumber
\end{eqnarray}

\begin{lemma}\label{atboundary}
{\em Consider a sequence $(\Lambda_k,\Theta_k), k\ge 1$ in ${\cal L}_+$ such that the matrix \\
$\lim_k\left(E_n\Lambda_k E_n\tp +{\Ub}_N \Theta_k {\Ub}_N \tp -\Theta_k\right)$ is singular. Assume also that the sequence $(\Lambda_k,\Theta_k)$ is bounded. Then, $J(\Lambda_k,\Theta_k)\rightarrow\infty$. }
\end{lemma}
\begin{IEEEproof}
Simply write 
$$J(\Lambda_k,\Theta_k)=-\log\det\left(E_n\Lambda_k E_n\tp +{\Ub}_N \Theta_k {\Ub}_N \tp -\Theta_k\right)+\Tr(\Lambda_k\Tb_k).$$ Since $\Tr(\Lambda_k\Tb_k)$ is bounded, the conclusion follows.
\end{IEEEproof}

\noindent
{\bf Proof of Theorem \ref{existencedual}.}
Observe that the function $J$ is a continuous, bounded below (Lemma \ref{boundedbelow}) function that tends to infinity both when  $\|(\Lambda,\Theta)\|$ tends to infinity (Lemma \ref{atinfinity}) and when it tends to the boundary $\partial{\cal L}_+$ with $\|(\Lambda,\Theta)\|$ remaining bounded (Lemma \ref{atboundary}). It follows that $J$ is {\em inf-compact} on ${\cal L}_+$, namely it has compact sublevel sets. By Weierstrass' Theorem\footnote{A continuous function on a compact set always achieves its  maximum and minimum on that set.}, it admits at least one minimum point. Since $J$ is strictly convex, the minimum point is unique.
\hfill $\Box$ \vskip 2ex

\section{Reconciliation with Dempster's Covariance Selection}\label{sec:Reconciliation}

Let $(\bar{\Lambda},\bar{\Theta})$ be the unique minimum point of $J$ in ${\cal L}_+$ (Theorem \ref{existencedual}). Then $\Sigmab_N^o\in{\mathfrak S}_{N,+}$ given by
\begin{equation}\label{optimumprimal}
\Sigmab_N^o=\left(E_n\bar{\Lambda} E_n\tp +{\Ub}_N \bar{\Theta} {\Ub}_N \tp -\bar{\Theta}\right)^{-1}
\end{equation}
satisfies (\ref{c1}) and (\ref{c2}). Hence, it is the unique solution of the primal problem (MEP). Since it satisfies (\ref{c2}), $\Sigmab_N^o$ is in particular a block-circulant matrix and hence  so is 
 $$
 (\Sigmab_N^o)^{-1}=\left(E_n\bar{\Lambda} E_n\tp +{\Ub}_N \bar{\Theta} {\Ub}_N \tp -\bar{\Theta}\right).
 $$
 Let $\pi_{\mathfrak{C}_N}$ denote the orthogonal projection onto the linear subspace of symmetric, block-circulant matrices $\mathfrak{C}_N$. It follows that, in force of Lemma \ref{perpcirc},
\begin{equation}\label{representation}(\Sigmab_N^o)^{-1}=\pi_{\mathfrak{C}_N}((\Sigmab_N^o)^{-1})=\pi_{\mathfrak{C}_N}\left(E_n\bar{\Lambda} E_n\tp +{\Ub}_N \bar{\Theta} {\Ub}_N \tp -\bar{\Theta}\right)=\pi_{\mathfrak{C}_N}\left(E_n\bar{\Lambda} E_n\tp \right)\,.
\end{equation}

\begin{theorem}\label{tridiagonalinverse}
{\em Let $\Sigmab_N^o$ be the maximum Gaussian entropy covariance given by (\ref{optimumprimal}). Then $(\Sigmab_N^o)^{-1}$ is a symmetric   block-circulant matrix which is banded of bandwidth $n$. Hence  the solution of (MEP) may be viewed as the covariance of a stationary reciprocal process of order $n$ defined on $\Zbb_N$.}
\end{theorem}
\begin{IEEEproof}
Let
$$
\Pi_{\bar{\Lambda}}:=\pi_{\mathfrak{C}_N}\left(E_n\bar{\Lambda} E_n\tp \right)=\bmat
\Pi_{0}&\Pi_{1}\tp &\Pi_{2}\tp &\dots&\Pi_{1}\\
\Pi_{1}&\Pi_{0}&\Pi_{1}\tp &\dots&\Pi_{2}\\
\vdots&\ddots&\ddots
&\ddots&\vdots\\
       \Pi_{2}\tp &\dots&\Pi_{1}&\Pi_{0}&\Pi_{1}\tp \\
       \Pi_{1}\tp &\Pi_{2}\tp &\dots& \Pi_{1}&\Pi_{0}
       \emat
$$
be the orthogonal projection of $\left(E_n\bar{\Lambda} E_n\tp \right)$ onto $\mathfrak{C}_N$. Since $\Pi_{\bar{\Lambda}}$ is symmetric and block-circulant, it is characterized by the orthogonality condition
\begin{equation}\label{orthogonality}
\Tr\left[\left(E_n\bar{\Lambda} E_n\tp -\Pi_{\bar{\Lambda}}\right)C\right]=\langle E_n\bar{\Lambda} E_n\tp -\Pi_{\bar{\Lambda}},C\rangle=0,\quad \forall C\in\mathfrak{C}_N.
\end{equation}
Next observe that, if we write $C=\Circ\left[C_{0},C_{1} ,C_{2} ,\ldots,C_{2}\tp,C_{1}\tp \right]$ and
$$
\bar{\Lambda}=\bmat 
\bar{\Lambda}_{00} &\bar{\Lambda}_{01}& \ldots&\ldots& \bar{\Lambda}_{0n}\\
\bar{\Lambda}_{10}\tp &\bar{\Lambda}_{11} &\ldots & &\bar{\Lambda}_{1n}\\
           \ldots  &             &        \ldots        &         & \ldots\\
\bar{\Lambda}_{n0}\tp &\bar{\Lambda}_{n1}\tp &\ldots & &\bar{\Lambda}_{nn} 
\emat,\qquad \bar{\Lambda}_{k,j}= \bar{\Lambda}_{j,k}\tp
$$
 then
\begin{eqnarray*}
\Tr \left[ E_n\bar{\Lambda} E_n\tp C \right] &=& \Tr \left[ \bar{\Lambda} E_n\tp CE_n\right]=\Tr \left[(\bar{\Lambda}_{00}+\bar{\Lambda}_{11} +\ldots +\bar{\Lambda}_{nn})C_{0} \right.\\
&+& (\bar{\Lambda}_{01} + \bar{\Lambda}_{12} +\ldots + \bar{\Lambda}_{n-1,n})C_{1}+ \ldots + \bar{\Lambda}_{0n}  C_{n} \\
&+& \left.(\bar{\Lambda}_{10} + \bar{\Lambda}_{21} +\ldots,\bar{\Lambda}_{n,n-1})C_{1}\tp+ \ldots + \bar{\Lambda}_{n0}  C_{n}\tp \right]\,.
\end{eqnarray*}
On the other hand, recalling that the product of two block-circulant matrices is block-circulant, we have that $\Tr \left[\Pi_{\bar{\Lambda}}C\right]$ is simply $N$ times the trace of the first block row of $\Pi_{\bar{\Lambda}}$ times the first block column of $C$. We get
\begin{equation*}
\Tr\left[\Pi_{\bar{\Lambda}}C\right]=N\Tr\left[\Pi_{0}C_{0}+\Pi_{1}\tp C_{1}+\Pi_{2}\tp C_{2}+  \ldots+\Pi_{2}C_{2}\tp +\Pi_{1}C_{1}\tp \right].
\end{equation*}
Hence, the orthogonality condition (\ref{orthogonality}), reads
\begin{eqnarray}\nonumber
\Tr\left[\left(E_n\bar{\Lambda} E_n\tp -\Pi_{\bar{\Lambda}}\right)C\right] &=&\Tr\left[\left((\bar{\Lambda}_{00}+\bar{\Lambda}_{11} +\ldots +\bar{\Lambda}_{nn}) -N\Pi_{0}\right)C_{0}+\right.\nonumber\\
&+&\left( (\bar{\Lambda}_{01} + \bar{\Lambda}_{12} +\ldots + \bar{\Lambda}_{n-1,n})-N\Pi_{1}\tp \right)C_{1}\nonumber \\
&+&\left((\bar{\Lambda}_{10} + \bar{\Lambda}_{21} +\ldots,\bar{\Lambda}_{n,n-1})-N\Pi_{1}\right)C_{1} \tp\nonumber \\
&+& \left.\ldots (\bar{\Lambda}_{0n}-N\Pi_{1}\tp)C_{n} +  (\bar{\Lambda}_{n0}-N\Pi_{1} )C_{n}\tp)\right] \nonumber \\
&+ &N\Pi_{n+1}\tp C_{n+1}+N\Pi_{n+1}C_{n+1}\tp + N\Pi_{n+2}\tp C_{n+2}+N\Pi_{n+2}C_{n+2}\tp +  \ldots  =0. \nonumber
\end{eqnarray}
Since this must hold true forall $C\in\mathfrak{C}_N$, we conclude that
\begin{eqnarray*}
&&\Pi_{0}=\frac{1}{N}\,(\bar{\Lambda}_{00}+\bar{\Lambda}_{11} +\ldots +\bar{\Lambda}_{nn}),\\
&&\Pi_{1}=\frac{1}{N}(\bar{\Lambda}_{01} + \bar{\Lambda}_{12} +\ldots + \bar{\Lambda}_{n-1,n})\tp,\\
& & \ldots\\
&&\Pi_{n}=\frac{1}{N} \bar{\Lambda}_{0n}\tp \,,
\end{eqnarray*}
while from the last equation we get $\Pi_i=0$,  forall $i$ in the interval $\, n+1\leq i \leq N-n-1\,$. From this it is clear that the inverse of the covariance matrix solving the primal problem (MEP), namely $\Pi_{\bar{\Lambda}}=(\Sigmab_N^o)^{-1}$ has a circulant block-banded structure of bandwidth $n$.
\end{IEEEproof}

Since the beginning of Section V, we have been dealing only with Gaussian distributions in order to facilitate the comparison with Dempster's classical results. It is now time to show that the Gaussian assumption can be dispensed with, and our solution is indeed optimal in the larger family of (zero-mean) second-order distributions.
\begin{theorem}\label{GeneralEntropy}
{\em The Gaussian distribution with (zero mean and) covariance  $\Sigmab_N^o$ defined by  (\ref{optimumprimal})    maximizes  the entropy functional \eqref{DiffEntropy}   over the set of all (zero mean)  probability densities whose covariance matrix satisfies the boundary conditions \eqref{c1}, \eqref{c2}.}
\end{theorem}
\begin{IEEEproof}
Let $\mathfrak{C}_N(\Tb_n)$ be the set of (block-circulant) covariance matrices satisfying the boundary conditions \eqref{c1}, \eqref{c2} and let $p_{\Sigmab}$ be a probability density with zero mean and covariance $\Sigmab$. In particular, we shall denote by $g_{\Sigmab}$ the Gaussian density with zero mean and covariance $\Sigmab$. Now, by a famous theorem of Shannon  \cite{Shannon-48}, the probability distribution having maximum entropy in the class of all distribution with a fixed mean vector (which we take equal to zero) and variance matrix  $\Sigmab$, is the Gaussian distribution $g_{\Sigmab}$. Hence:
$$
\max_{\Sigmab \in \mathfrak{C}_N(\Tb_n)}\, \left\{ \max_{p_{\Sigmab}}\, \left[H(p_{\Sigmab})\right]\, \right\}= \max_{\Sigmab \in \mathfrak{C}_N(\Tb_n)}\, \left\{ \, H(g_{\Sigmab})\, \right\}
$$
where the  maximum in the right-hand side is attained by $g_{\Sigmab_N^o}$.  
\end{IEEEproof}

The above   can be interpreted as  a particular    covariance selection result  in the vein of Dempster's paper; compare in particular  \cite[Proposition a]{Dempster-72}.  In  fact the results of this section substantiate also the maximum  entropy principle of Dempster (Proposition \ref{prop:b}).  It is however important to note that none of our results    follows   as a particular case  from   Dempster's results, since  \cite{Dempster-72} deals with a very  unstructured setting.  In particular our main result (Theorem \ref{tridiagonalinverse})  that the solution, $\Sigmab_N^o$, to our primal problem (MEP) has  a block-circulant {\em banded}  inverse, is completely original. Its proof    uses in an essential way the characterization of the MEP solution provided by our variational analysis and   cleverly exploits the block-circulant structure.

Actually, our results, {\em together} with Dempster's, may be used to show that the maximum entropy distribution, subject only to  moment constraints 
(compatible with the circulant structure) on a block band and on the corners, is necessarily block-circulant, i.e. the underlying process is stationary\footnote{An alternative proof of this fact can be constructed based on the invariance properties of the entropy functional and its strict concavity. This has been recently accomplished  (in a more general framework) in \cite{Carli-Georgiou}.}. 
 
 Because of the equivalence of reciprocal AR modeling and the underlying process covariance having an inverse with a banded structure, explained in Section \ref{sec:circext},  we see that the Maximum Entropy principle  leads in fact to (reciprocal) AR models.  This makes contact with the ever-present problem in control an signal processing of (approximate) AR modeling from finite covariance data, whose solution dates back to the work of N. Levinson and P. Whittle.  That AR modeling from finite covariance data is actually equivalent to a positive   band extension problems for infinite Toeplitz matrices has been realized and  studied in the past decades by Dym, Gohberg and co-workers, see e.g. \cite{Dym-G-81}, \cite{Gohberg-G-K-94} as representative references of a very large literature. We should stress here that    band extension problems for infinite Toeplitz matrices are invariably attacked and solved by factorization techniques,   but circulant matrices do not fit in the ``banded algebra'' framework used in the literature. Also, one should note that  the maximum entropy property  is usually presented in the literature as a final embellishment of a solution  which was already obtained by factorization techniques. Here, for the   circulant band  extension problem, factorization techniques do not   work and  the maximum  entropy principle turns out  to be the key to the solution of the problem.
 
This fact,  together with Dempster's observation \cite[Proposition b]{Dempster-72}, may    be taken as   a proof  (although referred to a very specific case) of a very much quoted general principle   that maximum entropy distributions are distributions achieving maximum {\em simplicity of explanation} of the data.

Finally,   we anticipate that the results of this section lead  to an efficient iterative algorithm for the explicit solution of the MEP which is guaranteed to converge to a unique minimum. This solves  the variational problem and hence   the circulant band extension problem which subsumes  maximum likelihood identification of reciprocal processes. This algorithm, which will not be described here for reasons of space limitations, compares very favorably with   the best techniques available so far.

\section{Conclusions}
A new class of stationary reciprocal processes on a finite interval  has been introduced which are the acausal analog of autoregressive (AR) processes on the integer line.   Maximum likelihood identification of these AR-type reciprocal models is discussed.  The computation of the estimates  of the matrix parameters   of the model turns out to be  a particular instance of a {\em Covariance selection problem} of the kind studied by the statistician     A.P. Dempster   in the early seventies. In matrix terminology, the  covariance selection   for  stationary reciprocal models is equivalent  to a special {\em matrix band extension problem} for block-circulant matrices.   We have shown that this band  extension problem can    be  solved by maximizing an entropy functional.

  \section*{Acknowledgements}
We are very grateful to the anonymous reviewers whose detailed comments and suggestions have helped us to improve the manuscript.

\bibliographystyle{plain}
\bibliography{biblio_recip}

\end{document}